\documentclass[11pt, twoside]{amsart}
\usepackage{amsfonts}

\begin{document}
\newtheorem{Mthm}{Main Theorem\!\!\,}
\newtheorem{Thm}{Theorem}[section]
\newtheorem{Prop}[Thm]{Proposition}
\newtheorem{Lem}[Thm]{Lemma}
\newtheorem{Cor}[Thm]{Corollary}
\newtheorem{Def}[Thm]{Definition}
\newtheorem{Guess}[Thm]{Conjecture}
\newtheorem{Ex}[Thm]{Example}
\newtheorem{Rmk}{Remark\!\!}
\newtheorem{Not}{Notation}

\renewcommand{\theThm} {\thesection.\arabic{Thm}}
\renewcommand{\theProp}{\thesection.\arabic{Prop}}
\renewcommand{\theLem}{\thesection.\arabic{Lem}}
\renewcommand{\theCor}{\thesection.\arabic{Cor}}
\renewcommand{\theDef}{\thesection.\arabic{Def}}
\renewcommand{\theGuess}{\thesection.\arabic{Guess}}
\renewcommand{\theEx}{\thesection.\arabic{Ex}}
\renewcommand{\theRmk}{}
\renewcommand{\theMthm}{}
\renewcommand{\theNot}{}
\renewcommand{\thefootnote}{\fnsymbol{footnote}}


\title[Toledo invariants of . . . homology 3-spheres]{\textbf\underline{Toledo invariants of Higgs bundles on elliptic surfaces associated to base orbifolds of Seifert fibered homology
3-spheres}}
\author{Mike Krebs}
\date{}
\maketitle



\begin{abstract}
To each connected component in the space of semisimple
representations from the orbifold fundamental group of the base
orbifold of a Seifert fibered homology 3-sphere into the Lie group
$\rm{U}(2,1)$, we associate a real number called the ``orbifold
Toledo invariant.''  For each such orbifold, there exists an
elliptic surface over it, called a Dolgachev surface.  Using the
theory of Higgs bundles on these Dolgachev surfaces, we explicitly
compute all values taken on by the orbifold Toledo invariant.
\end{abstract}

\section*{Introduction}

In this paper, we investigate the space
$\mathcal{R}^+_{\rm{U}(2,1)}(O)=\frac{\rm{Hom}^+(\pi^{\rm{orb}}_1(\it{O}),\,\rm{U}(2,1))}{\rm{U}(2,1)}$
of semisimple representations from the orbifold fundamental group
of a certain 2-orbifold $O$ into the Lie group $\rm{U}(2,1)$,
modulo conjugation.  To each connected component in
$\mathcal{R}^+_{\rm{U}(2,1)}(O)$, we associate a number that we
call the ``orbifold Toledo invariant.'' Our main result (Thm.
\ref{main_1}) explicitly computes all values that the orbifold
Toledo invariant takes on. One thereby obtains (Cor.
\ref{main_2}a) a lower bound for the number of connected
components in $\mathcal{R}^+_{\rm{U}(2,1)}(O)$.  The orbifolds we
consider are quotients of certain 3-manifolds $Y$---namely,
Seifert fibered homology 3-spheres---by the action of $S^1$. Our
results also yield (Cor. \ref{main_2}b) a lower bound for the
number of connected components in
$\mathcal{R}^{^*}\!\!_{\rm{PU}(2,1)}(Y)=\frac{\rm{Hom}^*(\pi_1(\it{Y}),\,\rm{PU}(2,1))}{\rm{PU}(2,1)}$,
the space of irreducible representations from the fundamental
group of $Y$ into $\rm{PU}(2,1)$, modulo conjugation.

In \cite{Toledo79}, Toledo introduces an invariant $\tau$ for
representations of the fundamental group of a oriented 2-manifold
$M$ into ${\rm PU}(p,1)$.  This invariant can be viewed as a map
$\tau:{\rm Hom}(\pi_1({M}),{\rm PU}(p,1))\rightarrow \mathbb{R}$.
As discussed in \S\ref{Toledo}, the construction of the Toledo
invariant is quite general: one may replace $M$ by an arbitrary
topological space and ${\rm PU}(p,1)$ by any topological group
$G$. Moreover, representations which take on distinct Toledo
invariants necessarily lie in distinct components of the
corresponding representation space.  In the case where $M$ is a
compact Riemann surface of genus $g>1$, previously established
results about the space
$\mathcal{R}_G^+(M)=\frac{\rm{Hom}^+(\pi_1(\it{M}),\,G)}{G}$ of
semisimple representations of $\pi_1(M)$ into $G$, modulo
conjugation, include:

\begin{itemize}

\item[$\bullet$] The Toledo invariant gives a bijection between
the set of all $\tau\in\frac{2}{3}\mathbb{Z}$ with $|\tau|\leq
2g-2$ and the set of all connected components in
$\mathcal{R}_{\rm{PU}(2,1)}^+(M)$ \cite{GKL}, \cite{Xia}.

\item[$\bullet$] If $\tau$ is sufficiently large and $c$ is any
integer, then the subset of
$\mathcal{R}_{\rm{PU}(\it{p},\it{p})}^+(M)$ corresponding to
representations with Toledo invariant $\tau$ and Chern class $c$
is connected \cite{Markman-Xia}.

\item[$\bullet$] The Toledo invariant gives a bijection between
the set of even integers $\tau$ with $|\tau|\leq~2(g-1)$ and the
set of connected components in
$\mathcal{R}_{\rm{U}(\it{p},\rm{1})}^+(M)$ \cite{Xia03}.

\item[$\bullet$] The subset $\mathcal{R}(\tau,c)$ of
$\mathcal{R}_{\rm{PU}(\it{p},\it{q})}^+(M)$ corresponding to
representations with Toledo invariant $\tau$ and Chern class $c$
is non-empty if and only if $$\tau=\frac{|qa-p(c-a)|}{p+q}\leq
(g-1)\cdot{\rm min}\{p,q\}$$ for some integer $a$. Moreover, if
this inequality is satisfied and $p+q$ and $c$ are coprime, then
$\mathcal{R}(\tau,c)$ is connected \cite{BGPG}.

\end{itemize}

Other results concerning Toledo invariants can be found in
\cite{BGPG05}, \cite{BILW}, \cite{BIW}, \cite{Goldman},
\cite{GKL}, \cite{GP00}, \cite{GP03}, \cite{Toledo79}, and
\cite{Toledo89}.

To our orbifold $O$ we associate a complex surface $X$, called a
Dolgachev surface, whose fundamental group is isomophic to
$\pi_1^{\rm{orb}}(O)$.  The reason for doing so is that for
complex algebraic manifolds $M$, we have a correspondence between
representations of $\pi_1(M)$ and certain algebro-geometric
objects on $M$ called Higgs bundles.  (A Higgs bundle on $M$
consists of a holomorphic vector bundle plus some extra data; see
\S\ref{Higgs} for the definition and basic properties.)  The
relationship between representations of $\pi_1(M)$ and holomorphic
vector bundles on $M$ has been developed over the last forty years
by Narasimhan and Seshadri \cite{NS}, Atiyah and Bott \cite{AB},
Hitchin \cite{Hitchin}, Donaldson \cite{Donaldson}, Corlette
\cite{Corlette}, Simpson \cite{Simpson}, and others.

In \S\ref{Hodge}, we obtain detailed information about those Higgs
bundles on the Dolgachev surface $X$ that correspond to semisimple
representations $\rho:\pi_1(X)\to{\rm U}(2,1)$.  In \cite{Xia},
Xia computes the Toledo invariant of such a representation in
terms of the associated Higgs bundle.  This computation, together
with the results of \S\ref{Hodge}, enables us to determine all
Toledo invariants that arise from semisimple representations
$\rho$.  In \S\ref{orbifold}, we define ``orbifold Toledo
invariants'' for representations of the orbifold fundamental group
of the 2-orbifold $O$ into $\rm{PU}(2,1)$ and show that these are
in one-to-one correspondence with Toledo invariants on the
Dolgachev surface $X$.  In \S\ref{main}, we put the pieces
together, obtaining numerical conditions which completely
determine whether or not a real number $\tau$ represents an
orbifold Toledo invariant:

\begin{Mthm} Let $O$ be the base orbifold of Seifert fibered homology
$3$-sphere $Y$ such that $\pi_1^{\rm{orb}}(O)$ is infinite.  Let
$n$ equal the number of cone points that $O$ has, and let
$m_1,\dots,m_n$ denote the orders of these cone points.  Let
$\tau\in\mathbb{R}$. Then there exists a semisimple representation
$\rho:\pi_1^{\rm{orb}}(O)\rightarrow \rm{U(2,1)}$ such that $\tau$
is the orbifold Toledo invariant of $\rho$ if and only if
$\tau=\pm(y+\sum\frac{y_k}{m_k})$ for some integers $(y, y_1,
\dots, y_n)$ with $0\leq y_k<m_k$ such that at least one of the
numerical conditions {\rm (i)--(iv)} holds:

\vspace{.1in}

\begin{list}{-}
{\setlength{\labelsep}{.1in} \setlength{\labelwidth}{.3in}
\setlength{\leftmargin}{.4in} \setlength{\itemindent}{0in}
\setlength{\rightmargin}{.2in}}

\item[{\rm (i)}] There exist integers $a, a_1, \dots, a_n, b, b_1,
\dots, b_n$ with $0\leq a_k,b_k<m_k$ such that $b\leq -2$, and
$a+\#\{k\thinspace |\thinspace a_k\neq 0\}\geq 2$, and $2A<B$, and
$A<2B$, and {\rm ($\star$)} below holds. \\

\vspace{-.12in}

\item[{\rm (ii)}] There exist integers $a, a_1, \dots, a_n, b,
b_1, \dots, b_n$ with $0\leq a_k,b_k<m_k$ such that
$-B<A<\frac{1}{2}B$, and $d_2\leq -2$, and $b\leq -2$,
and {\rm ($\star$)} below holds, and $d_1+1\leq {\rm min}(-d_2-1,-d_3-1)$ for every $(n+1)$-tuple of integers $(c,c_1,\dots,c_n)$ such that $0\leq c_k<m_k$ for all $k$ and $d_1\geq 0$ and $C\geq\frac{2}{3}(A+B)$. \\

\vspace{-.1in}

\item[{\rm (iii)}] $y+\sum\frac{y_k}{m_k}>0$, and
$\;2y+\#\{y_k\geq\frac{m_k}{2}\}\leq
-2$.\\

\vspace{-.1in}

\item[{\rm (iv)}] $y=y_k=0$ for all $k$. \\

\vspace{-.1in}

\item[{\rm ($\star$)}]
$3y+\sum\lfloor\frac{3y_k}{m_k}\rfloor=a+b$, and
$3y_k-\lfloor\frac{3y_k}{m_k}\rfloor m_k=a_k+b_k$ for
$k=1,\dots,n$.

\end{list}

Here we have used the notations $A=a+\sum\frac{a_k}{m_k}$;
$B=b+\sum\frac{b_k}{m_k}$; $C=c+\sum\frac{c_k}{m_k}$;
$d_1=b-c-\#\{b_k<c_k\}$; $d_2=a-b-\#\{a_k<b_k\}$; and
$d_3=a-c-\#\{a_k<c_k\}$.

\label{main_1_intro}
\end{Mthm}

As a corollary, we obtain a lower bound for the number of
connected components in $\mathcal{R}^+_{\rm{U}(2,1)}(O)$.  In
\S\ref{Seifert}, we show that irreducible $\rm{PU}(2,1)$
representations of $\pi_1(Y)$ are in one-to-one correspondence
with irreducible $\rm{PU}(2,1)$ representations of
$\pi^{\rm{orb}}_1(O)$.  The main theorem therefore also furnishes
a lower bound for the number of connected components in
$\mathcal{R}^{^*}\!\!_{\rm{PU}(2,1)}(Y)$.

The space of irreducible $\rm{SU}(2)$ representations of
$\pi_1(Y)$ has been studied in detail by Fintushel and Stern
\cite{Finster}, Bauer and Okonek \cite{BO}, Kirk and Klassen
\cite{KK}, Furuta and Steer \cite{FS}, Bauer \cite{Bauer}, and
Boden \cite{Boden}. (The motivation of these authors was the study
of the $\rm{SU}(2)$ Casson's invariant and Floer homology for such
spaces $Y$.)  In many of these papers, the method is to associate
to $Y$ an auxiliary object whose fundamental group is closely
related to that of $Y$.  In \cite{FS} and \cite{Boden}, the
auxiliary object is a 2-orbifold; they study parabolic Higgs
bundles on the orbifold's underlying Riemann surface, i.e.
$\mathbb{CP}^1$.  In this paper, as in \cite{BO}, the auxiliary
object is a Dolgachev surface.  (Along the same lines, in
\cite{Biswas}, an elliptic surface is used to study vector bundles
on a 2-orbifold; it is not clear, however, that the associated
elliptic surface is algebraic, as claimed---see, for example,
\cite[Example 13.2]{BPV}.)

One motivation for studying $\rm{PU}(2,1)$ representations of the
fundamental groups of 3-manifolds comes from spherical CR
geometry. A spherical CR structure on a 3-manifold $M$ is a system
of coordinate charts into $S^3$ so that the transition functions
are elements of ${\rm PU}(2,1)$.  (Here we regard $\rm{PU}(2,1)$
as the isometry group of the complex ball in $\mathbb{C}^2$ and
the conformal group of its boundary $S^3$.  See \cite{Goldman}.)
In \cite{KT}, Kamishima and Tsuboi classify those closed
orientable 3-manifolds that admit $S^1$-invariant spherical CR
structures; these include the Seifert-fibered homology 3-spheres
considered here.  The space $\frac{{\rm Hom}(\pi_1(M),{\rm
PU}(2,1))}{{\rm PU}(2,1)}$ provides a local model for the
deformation space of spherical CR structures on $M$ \cite{KTan}.

Our lower bound for the number of components in
$\mathcal{R}^{^*}\!\!_{\rm{PU}(2,1)}(Y)$ takes into account only
those ${\rm PU}(2,1)$ representations which lift to ${\rm U}(2,1)$
representations.  Moreover, for $\mathcal{R}^+_{\rm{U}(2,1)}(O)$,
we conjecture that the number of components is in general strictly
greater than the number of orbifold Toledo invariants that occur.
We plan to continue investigating these representation spaces,
with the goal of precisely determining the number of components in
them.

The author is grateful to I. Dolgachev, E. Falbel, H. Ren, R.
Seyyedali, S. Zrebiec, G. Tinaglia, and most especially his thesis
advisor, Richard Wentworth, for many helpful discussions.

\section{Toledo invariants} \label{Toledo}

Given a manifold (or topological space) $M$ and a topological
group $G$, one may wish to study the representation variety
$\mathcal{R}=\frac{\rm{Hom}(\pi_1(\it{M}),G)}{G}$.  The goal of
this section is to define a family of invariants, called Toledo
invariants, that can be used to distinguish components of
$\mathcal{R}$. We then describe one such Toledo invariant more
specifically in the case where $G=U(2,1)$.

\subsection{The ``abstract nonsense'' of Toledo invariants}

Let $B$ be a solid topological space \cite{Steenrod}.  (Euclidean
space $\mathbb{R}^n$ is solid, for example.)  Let $G$ be a
topological group acting continuously on $B$ on the left.  We now
take $\omega$ to be a fixed $G$-invariant representative of a
cohomology class in $H^{^*}\!(B,\mathbb{C})$.  (If $B$ is a
manifold, we may regard $\omega$ as a closed singular cochain or
as a closed differential form, depending which is more
convenient.)

Let $M$ be a $C^{\infty}$ manifold.  We define a map
$\tau^{B,G,\omega}$ from Hom$(\pi_1(M),G)$ to
$H^{^*}\!(M,\mathbb{C})$ as follows.  Let $\rho\in$
Hom$(\pi_1(M),G)$. Let $\tilde M$ be the universal cover of $M$.
Note that $\pi_1(M)$ acts on $\tilde M\times B$ by
$\gamma\cdot(m,x)=(\gamma\cdot m,\rho(\gamma)\cdot x)$. Let
$E_{\rho}$ be the flat $B$-bundle on $M$ obtained by taking
$\tilde M\times B$ modulo the action of $\pi_1(M)$.  Let
$\pi_B:\tilde M\times B\rightarrow B$ be the projection map onto
the second factor, and let $\varphi$ be the natural map from
$\tilde M\times B$ to $E_{\rho}$. Since $\pi_1(M)$ acts freely on
$\tilde M$ and $\omega$ is $G$-invariant and closed, the pullback
$\pi_B^*\omega$ descends to $E_{\rho}$, where it represents a
cohomology class $[\varphi_*\pi_B^*\omega]\in
H^{^*}\!(E_{\rho},\mathbb{C})$.  Since the fibre $B$ is solid,
$E_{\rho}$ has a section; moreover, any two sections are homotopic
\cite[Theorem 12.2]{Steenrod}. Consequently,
$[s^*\varphi_*\pi_B^*\omega]$ is a well-defined cohomology class
in $H^{^*}\!(M,\mathbb{C})$.

\begin{Def}
The \emph{Toledo invariant} $\tau^{B,G,\omega}(\rho)$ is defined
by $$\tau^{B,G,\omega}(\rho)=[s^*\varphi_*\pi_B^*\omega]$$ for any
section $s$ of $E_{\rho}$.

\label{Toledo_0.5}
\end{Def}

\begin{Lem}
Let $M$ be a $C^{\infty}$ manifold, let
$\rho\in\,$\emph{Hom}$(\pi_1(M),G)$, let $g\in G$, and define
$\rho\, ':\pi_1(M)\rightarrow\pi_1(M)$ by $\rho\,
'(\gamma)=g\rho(\gamma)g^{-1}$.  Then
$\tau^{B,G,\omega}(\rho)=\tau^{B,G,\omega}(\rho\, ')$.  In other
words, the Toledo invariant is invariant under conjugation.

\label{Toledo_1}
\end{Lem}
\emph{Proof.}  We define a map $\psi:\tilde M\times B\rightarrow
\tilde M\times B$ by $\psi(x,b)=(x,g\cdot b)$.  Let
$E_{\rho}=\frac{\tilde M\times B}{\pi_1(M)}$ (where the action is
induced by $\rho$), and let $E_{\rho\, '}=\frac{\tilde M\times
B}{\pi_1(M)}$ (where the action is induced by $\rho\, '$).  Then
$\psi$ descends to a map from $E_{\rho}$ to $E_{\rho\, '}$; we
denote this new map by $\psi$ as well.  If $s$ is a section of
$E_{\rho}$, then $s'=\psi\circ s$ is a section of $E_{\rho\, '}$.
In summary, we have that the following diagram commutes.

\vspace{.2in}

\begin{picture}(360,180)

\put(110,0){$M$}

\put(252,0){$E_{\rho\, '}$}

\put(0,80){$M$}

\put(135,80){$E_{\rho}$}

\put(240,120){$\tilde M\times B$}

\put(356,120){$B$}

\put(120,175){$\tilde M\times B$}

\put(245,175){$B$}

\put(15,70){\vector(3,-2){90}}

\put(50,30){$id$}

\put(130,5){\vector(1,0){110}}

\put(175,10){$s'$}

\put(155,75){\vector(3,-2){90}}

\put(190,57){$\psi$}

\put(20,85){\vector(1,0){103}}

\put(65,90){$s$}

\put(140,170){\vector(0,-1){73}}

\put(125,135){$\varphi$}

\put(260,110){\vector(0,-1){90}}

\put(270,70){$\varphi'$}

\put(150,170){\vector(2,-1){80}}

\put(185,137){$\psi$}

\put(260,175){\vector(2,-1){90}}

\put(307,157){$id$}

\put(160,179){\vector(1,0){80}}

\put(200,170){$\pi_B$}

\put(280,123){\vector(1,0){70}}

\put(310,114){$\pi_B$}

\end{picture}

\vspace*{.1in}

The lemma follows from chasing this diagram.$\;\square$

\vspace*{.1in}

Let $G$ act on $\mbox{Hom}(\pi_1(M),G)$ by conjugation.  Lemma
\ref{Toledo_1} shows that $\tau^{B,G,\omega}(\rho)$ defines a map
from $\frac{{\rm Hom}(\pi_1(M),G)}{G}$ to
$H^{^*}\!(M,\mathbb{C})$. Topologize $\frac{{\rm
Hom}(\pi_1(M),G)}{G}$ by giving $\mbox{Hom}(\pi_1(M),G)$ by the
point-open topology and giving $\frac{{\rm Hom}(\pi_1(M),G)}{G}$
the quotient topology. Note that if $t_1,\dots,t_n$ are generators
for $\pi_1(M)$, then ${\rm Hom}(\pi_1(M),G)$ is homeomorphic to
the closed subspace $\{(x_1,\dots,x_n)\in
G^n\;|\;r_{\alpha}(x_1,\dots,x_n)=1\}$ of $G^n$, where the
$r_{\alpha}$'s range over all relations between the $t$'s.

\begin{Lem}
Suppose that $B$ and $M$ are $C^{\infty}$ manifolds, that $M$ is
compact, that $G$ is a Lie group, and that $\omega$ is a closed
$G$-invariant $k$-form on $B$. Topologize $H^k(M,\mathbb{C})$ as a
finite-dimensional vector space.  Then $\tau^{B,G,\omega}$ defines
a continuous function from $\frac{{\rm Hom}(\pi_1(M),G)}{G}$ to
$H^k(M,\mathbb{C})$.

\label{Toledo_2}
\end{Lem}
\emph{Proof.}  It suffices to show that $\tau^{B,G,\omega}$ is
continuous on $\mbox{Hom}(\pi_1(M),G)$.  Let $C={\rm
Hom}(\pi_1(M),G)\times\tilde M\times B$.  An action of $\pi_1(M)$
on $C$ is given by
$$\gamma\cdot(\rho,m,x)=(\rho,\gamma\cdot~m,\rho(\gamma)\cdot
x).$$ Then $\frac{C}{\pi_1{M}}$ is a fibre bundle over
$\mbox{Hom}(\pi_1(M),G)\times M$ with fibre $B$. Since $B$ is
solid and since ${\rm Hom}(\pi_1(M),G)$ is a subspace of $G^n$,
there exists a section $s:\mbox{Hom}(\pi_1(M),G)\times
M\rightarrow\frac{C}{\pi_1{M}}$. Lift $s$ to a map $\tilde s:{\rm
Hom}(\pi_1(M),G)\times\tilde M\to C$.  By Tietze extension and by
inclusion of $C$ into $G^n\times\tilde M\times B$, we have that
$\tilde s$ extends to a map $\tilde s$ from $G^n\times\tilde M$ to
$G^n\times\tilde M\times B$.  Let $\pi_B:G^n\times\tilde M\times
B\to B$ denote projection onto the third factor.  Given $\rho\in
G^n$, define $\iota_{\rho}:\tilde M\rightarrow G^n\times M$ by
$\iota_{\rho}(m)=(\rho,m)$.  Then $\iota_{\rho}^*\tilde
s^*\pi_B^*\omega$ defines a closed $\pi_1(M)$-invariant singular
$k$-cochain on $\tilde M$, which therefore defines a closed
cochain $\tau_0(\rho)$ on $M$.  Let $\tau(\rho)$ be the associated
cohomology class in $H^k(M,\mathbb{C})$. Note that $\tilde
s\circ\iota_{\rho}$ defines a section of $E_{\rho}$, the fibre
bundle from the definition of the Toledo invariant. Therefore, we
have that $\tau^{B,G,\omega}(\rho)=\tau(\rho)$.  We now show that
$\tau(\rho)$ varies continuously with $\rho$.

\vspace*{.1in}

Let $r$ be the dimension over $\mathbb{C}$ of $H^k(M,\mathbb{C})$.
Given open sets $U_1,\dots,U_r$ of $\mathbb{C}$ and closed
cochains $\sigma_1,\dots,\sigma_r\in C^k(M,\mathbb{C})$ such that
the associated cohomology classes $[\sigma_1],\dots,[\sigma_r]$
are linearly independent, let $IB(\sigma_1,\dots,\sigma_r,
U_1,\dots,U_r)=\{\sum a_{\ell}\sigma_{\ell}\;|\;a_{\ell}\in
U_{\ell}\}.$  Then $\tau$ is ${\rm continuous}$ if and only if
$\tau_0^{-1}(IB(\sigma_1,\dots,\sigma_r, U_1,\dots,U_r))$ is open
for all\linebreak $\sigma_1, \dots, \sigma_r, U_1,\dots,U_r$.

\vspace*{.1in}

Now fix $r$ closed cochains $\sigma_1,\dots,\sigma_r\in
C^k(M,\mathbb{C})$ such that the associated cohomology classes
$[\sigma_1],\dots,[\sigma_r]$ are linearly independent.  Let
$C_1,\dots C_r$ be singular $k$-simplices in $M$ such that the map
from $\phi:IB(\sigma_1,\dots,\sigma_r,
U_1,\dots,U_r)\to\mathbb{C}^r$ defined by
$\phi(\sigma)=(<\sigma,C_1>,\dots,<\sigma,C_r>)$ is bijective.
Define $\phi_{\ell}:IB(\sigma_1,\dots,\sigma_r,
U_1,\dots,U_r)\to\mathbb{C}$ by
$\phi_{\ell}(\sigma)=<\!\sigma,C_{\ell}\!>$. It now suffices to
show that $\tau_0^{-1}(\phi_{\ell}^{-1})(U_{\ell})$ is open for
every open set $U_{\ell}$ of $\mathbb{C}$.

\vspace*{.1in}

By subdividing $C_{\ell}$ into small enough pieces, we can assume
that the image of $C_{\ell}$ is a subset of an open set $V$ of $M$
such that $V$ is homeomorphic, via the natural covering map, to an
open set $\tilde V$ of $\tilde M$.  We may then regard $C_{\ell}$
as a map from the standard $k$-simplex $\Delta_k$ to $\tilde M$.

\vspace*{.1in}

Let $\rho_0\in\tau_0^{-1}(\phi_{\ell}^{-1})(U_{\ell})$.  Endow
$G^n$ with a Riemannian metric.  We now show that for sufficiently
small $\delta$, the ball of radius $\delta$ centered at $\rho_0$
lies entirely within $\tau_0^{-1}(\phi_{\ell}^{-1})(U_{\ell})$;
this will conclude our proof.  Let $\rho_1\in G^n$. Let $c(t)$ be
a geodesic in $G^n$ with $c(0)=\rho_0$ and $c(1)=\rho_1$.  Let
$h:\Delta_k\times [0,1]\to B$ be a piecewise smooth function
homotopic to $\pi_B\circ\tilde s\circ
\iota_{c(t)}\circ(C_{\ell}\times{\rm id})$.  Note that
$\phi_{\ell}(\tau(\rho_j))=\int_{\Delta_k\times \{j\}} h^*\omega$
for $j=0,1$.  Stokes' Theorem then guarantees, for any
$\epsilon>0$, the existence of a $\delta>0$ such that
$|\phi_{\ell}(\tau(\rho_1))-\phi_{\ell}(\tau(\rho_0))|<\epsilon$
if the distance from $\rho_0$ to $\rho_1$ is less than
$\delta$.$\;\;\;\square$

\vspace*{.1in}

\textbf{Remark.}  If the image of $\tau^{B,G,\omega}$ is discrete,
then Lemma \ref{Toledo_2} shows that $\tau^{B,G,\omega}$ is
constant on connected components of $\frac{{\rm
Hom}(\pi_1(M),G)}{G}$.  This will be the case in our main theorem
(Thm. \ref{main_1}); the number of distinct values of
$\tau^{B,G,\omega}$ therefore provides, in this case, a lower
bound for the number of connected components in $\frac{{\rm
Hom}(\pi_1(M),G)}{G}$.  Lemma \ref{Toledo_2} is used (implicitly)
in this manner in \cite{Xia}, \cite{Xia03}, \cite{Markman-Xia},
and \cite{BGPG}.

\vspace*{.1in}

\textbf{Example.}  A simple example shows that $\tau^{B,G,\omega}$
is not always constant on connected components of ${\rm
Hom}(\pi_1(M),G)$.  Let $M$ be the unit circle $S^1$, let
$G=B=\mathbb{R}$ (where $G$ acts on $B$ by translation), and let
$\omega=dx$.  Let $t$ be the standard generator of $\pi_1(M)$, and
identify ${\rm Hom}(\pi_1(M),G)$ with $\mathbb{R}$ by
$\rho\mapsto\rho(t)$. Since ${\rm Hom}(\pi_1(M),G)$ has a single
connected component, it suffices to show that the Toledo invariant
is not a constant function.  Identifying $\tilde M$ with
$\mathbb{R}$ in the usual way, a $\rho$-equivariant section of
$\tilde M\times B$ is given by $x\mapsto(x,\rho(t)x)$.  One can
then compute that the Toledo invariant $\tau^{B,G,\omega}(\rho)$
is the cohomology class defined by
$\rho(t)d\theta$.$\;\;\;\square$

\subsection{The ${\rm U}(2,1)$ Toledo invariant}

We now turn our attention to the special case of this construction
that will be the focus of the remainder of this paper. Define
$g:\mathbb{C}^3\rightarrow\mathbb{C}$ by
$g(z_0,z_1,z_2)=|z_0|^2-|z_1|^2-|z_2|^2.$  Let ${\rm U}(2,1)=\{
A\in {\rm GL}(3,\mathbb{C})\,|\;g(Az)=g(z)\;\mbox{for
all}\;z\in\mathbb{C}^3\}$.

\vspace*{.1in}

Let $G=\mbox{U}(2,1)$.  Define $B$ by:
$$B=\textbf{H}_{\mathbb{C}}^2=\{(1,z_1,z_2)\in\mathbb{C}^3:1-|z_1|^2-|z_2|^2<1\}.$$
($B$ is the ball model of 2-dimensional complex hyperbolic space
\cite{Goldman}.)  Note that $B$ is homeomorphic to $\mathbb{R}^4$,
hence solid.  $G$ acts on $B$ as follows.  Let $z\in B$ and $A\in
G$.  Define the action of $A$ on $z$ by $A\cdot
z=\lambda\cdot(Az)$, where the $Az$ on the right hand side is
given by ordinary matrix multiplication (regarding $z$ as a column
vector), and $\lambda$ is the unique complex number such that
$\lambda\cdot(Az)\in B$. (We know that $\lambda$ exists since
U(2,1) preserves the indefinite form $|z_0|^2-|z_1|^2-|z_2|^2$.)
Let $\omega=\frac{i}{2\pi}\partial\overline\partial\log g$. Note
that $\omega$ is invariant under multiplication by elements of
U(2,1) (since $g$ is) and is invariant under multiplication by
scalars. By the definition of the action of $G$ on $B$, then, the
restriction of $\omega$ to $B$ is $G$-invariant.  The center
$Z(\rm{U}(2,1))$ of $\rm{U}(2,1)$ equals $\{\lambda
I|\lambda\in\rm{U}(1)\}$. Let
$\rm{PU}(2,1)=\frac{\rm{U}(2,1)}{\it{Z}(\rm{U}(2,1))}$.  Then
there is an action of $\rm{PU}(2,1)$ on $B$, inherited from the
$\rm{U}(2,1)$ action.  It follows that $\omega$ is
$\rm{PU}(2,1)$-invariant.  From now on, all Toledo invariants will
have $B$ and $\omega$ as in this paragraph and $G=\rm{U}(2,1)$ or
$G=\rm{PU}(2,1)$.

\section{PU(2,1) representations of fundamental groups of Seifert fibered homology 3-spheres}
\label{Seifert}

The goal of this section is to note the relationship between
PU(2,1) representations of the fundamental group of a Seifert
fibered homology 3-sphere and PU(2,1) representations of the
fundamental group of a certain elliptic surface called a Dolgachev
surface.

Let $Y$ be a Seifert fibered homology 3-sphere.  (Following Lemma
\ref{Seifert_0.5}, we shall impose some additional constraints on
$Y$.). For the definition of Seifert fibered spaces and basic
facts about them, we refer to \cite{Orlik}.  A $(2n+1)$-tuple
$(-c_0; (m_1,c_1),\dots,(m_n,c_n))$ of integers, with $m_k$
positive for all $k$, is associated to $Y$. (These integers are
called the Seifert invariants of $Y$; we may think of $m_k$ as the
degree of twisting of the $k$th singular fibre of $Y$.)  To be a
homology 3-sphere, we must have that ${\rm gcd}(m_j,m_k)=1$
whenever $j\neq k$ \cite{FS}.  The notations $Y,n,$ and $(-c_0;
(m_1,c_1),\dots,(m_n,c_n))$ will be fixed throughout the rest of
this paper.

The fundamental group of $Y$ has the following presentation
\cite[section 5.3]{Orlik}:
\[\pi_1(Y)=\,<t_1,\dots,t_n,h\,|\, t_k^{m_k}h^{c_k}=t_1\dots t_n
h^{c_0}=[h,t_k]=1>\]

If $G$ is any group, then let $Z(G)$ denote its center.  We have
that $Z(\pi_1(Y))$ is generated by $h$ \cite[section 5.3]{Orlik},
so

\[\frac{\pi_{1}(Y)}{Z(\pi_{1}(Y))}=\,<t_1,\dots,t_n\,|\, t_k^{m_k}=t_1\dots t_n=1>.\]

We now construct a complex surface $X$, called a Dolgachev
surface. The following description of this construction is taken
from \cite{BO}.  A generic cubic pencil in $\mathbb{CP}^2$ has
nine base points. Blowing up at these nine points, we obtain an
algebraic surface $X_0$ along with an elliptic fibration
$\pi_0:X_0\rightarrow\mathbb{CP}^1$.  Apply logarithmic
transformations \cite{GH} along $n$ disjoint nonsingular fibres of
$X_0$ with multiplicities $m_1,\dots,m_n$. The result is an
elliptic fibration $\pi:X\rightarrow\mathbb{CP}^{1}$, where $X$ is
the desired complex surface.  Throughout this paper, $X$ will
denote a Dolgachev surface whose invariants are $(m_1,\dots,m_n)$.

\begin{Lem}
$\pi_{1}(X)=\frac{\pi_{1}(Y)}{Z(\pi_{1}(Y))}.$  If $n\leq 2$, then
$\pi_1(X)$ is trivial.  If $n=3$ and $\{m_1,m_2,m_3\}=\{2,3,5\}$,
then $\pi_1(X)$ is the alternating group $A_5$.

\label{Seifert_0.5}
\end{Lem}
\emph{Proof.} \cite[Chapter II, \S 3]{Dolgachev} or \cite[Prop.
1.2 and subsequent discussion]{BO}.$\;\;\;\;\square$

Because of Lemma \ref{Seifert_0.5}, we will impose the
restrictions that $n\geq 3$ and that if $n=3$, then
$\{m_1,m_2,m_3\}\neq\{2,3,5\}$.

\begin{Def}
The Lie group \emph{PU(2,1)} acts on its Lie algebra $\frak{g}$
via the adjoint representation.  Consequently, if $H$ is a group
and $\rho:H\rightarrow \emph{PU}(2,1)$ is a representation, then
$\rho$ induces an action of $H$ on $\frak{g}$.  We say that $\rho$
is \emph{irreducible} (resp. \emph{reducible}) if this induced
action is irreducible (resp. reducible). We denote the set of
irreducible representations $\rho:H\rightarrow \emph{PU}(2,1)$ by
$\emph{Hom}^*(H,\emph{PU}(2,1))$.\end{Def}

\begin{Lem}
Let $H$ be a group, and let
$\rho\in\emph{Hom}^*(H,\emph{PU}(2,1))$. Then no points and no
complex geodesics in $\emph{\textbf{H}}_{\mathbb{C}}^2$ are
invariant under the action of $H$ on
$\emph{\textbf{H}}_{\mathbb{C}}^2$ induced by $\rho$.

\label{Seifert_1}
\end{Lem}
\emph{Proof.}  First, suppose that there exists
$x\in\textbf{H}_{\mathbb{C}}^2$ such that $\rho(h)\cdot x=x$ for
all $h\in H$.  Let $K=\{\phi\in\mbox{PU}(2,1)|\phi(x)=x\}$.  Then
$K$ is a Lie subgroup of PU(2,1); in fact, $K$ is conjugate to
P(U(2)$\times$U(1)).  Let $\frak{g}$ be the Lie algebra of
PU(2,1), and let $\frak{k}$ be the Lie subalgebra of $\frak{g}$
corresponding to $K$.  Since $\rho(H)\subset K$, we have that
$\frak{k}$ is invariant under the action of $H$ on
$\frak{g}$---but this is a contradiction, since $\rho$ is
irreducible.

\vspace*{.1in}

Similarly, suppose that $P$ is a complex geodesic in
$\textbf{H}_{\mathbb{C}}^2$ such that $\rho(h)\cdot x\in P$ for
all $h\in H$ and $x\in P$.  In this case, we take $K$ to be the
set of all elements in PU(2,1) that preserve $P$.  Again, $K$ is a
Lie subgroup of PU(2,1); this time, $K$ is conjugate to
P(U(1)$\times$U(1,1)).  Again, we find that $\frak{k}$ is
invariant under $H$, contradicting $\rho$'s
irreducibility.$\;\square$

\begin{Lem}
There exists a bijection $\varphi:{\rm Hom}^*(\pi_1(Y),{\rm
PU}(2,1))\rightarrow{\rm Hom}^*(\pi_1(X),{\rm PU}(2,1))$.

\label{Seifert_2}
\end{Lem}
\emph{Proof.}  Since
$\pi_{1}(X)=\frac{\pi_{1}(Y)}{Z(\pi_{1}(Y))}$, we have a
surjection $\sigma:\pi_1(Y)\rightarrow\pi_1(X)$, which in turn
induces an injection $\overline{\varphi}:{\rm Hom}(\pi_1(X),{\rm
PU}(2,1))\rightarrow{\rm Hom}(\pi_1(Y),{\rm PU}(2,1))$. Now,
$\rho$ and $\overline{\varphi}(\rho)=\sigma\circ\rho$ have the
same image, so $\rho$ is irreducible if and only if
$\overline{\varphi}(\rho)$ is irreducible.  Restricting
$\overline{\varphi}$ to the irreducible representations then gives
us an injection $\varphi$ from ${\rm Hom}^*(\pi_1(X),{\rm
PU}(2,1))$ to ${\rm Hom}^*(\pi_1(Y),{\rm PU}(2,1))$.  We must now
show that $\varphi$ surjects onto ${\rm Hom}^*(\pi_1(Y),{\rm
PU}(2,1))$.

Let $\tilde\rho:\pi_1(Y)\rightarrow{\rm PU}(2,1)$ be an
irreducible representation. We must find
$$\rho:\frac{\pi_{1}(Y)}{Z(\pi_{1}(Y))}\rightarrow{\rm PU}(2,1)$$
such that $\tilde\rho=\sigma\circ\rho$.  Recalling that the center
of $\pi_1(Y)$ is generated by the single element $h$, we see that
it suffices to prove that $\tilde\rho$ maps $h$ to the identity
element in ${\rm PU}(2,1)$.

Regard ${\rm PU}(2,1)$ as the group of isometries of
$\textbf{H}_{\mathbb{C}}^2$.  Our first goal is to show that
$\tilde\rho(h)$ has three linearly independent fixed points
$x_1,x_2,x_3$. Goldman \cite[p. 203]{Goldman} shows that
$\tilde\rho(h)$ has a fixed point
$x_1\in\textbf{H}_{\mathbb{C}}^2\cup\partial\textbf{H}_{\mathbb{C}}^2$.
 There must exist $f\in\tilde\rho(\pi_1(Y))$ such that $x_2=f(x_1)\neq
x_1$, else $\tilde\rho$ would not be irreducible, by Lemma
\ref{Seifert_1}. Since $h$ is central, $\tilde\rho(h)$ commutes
with $f$. Thus, $\tilde\rho(h)(x_2)=f(\tilde\rho(h)(x_1))=x_2.$
That is, $x_2$ is another fixed point of $\tilde\rho(h)$.  Let $P$
be the complex geodesic spanned by $x_1$ and $x_2$.  By linearity,
$P$ is invariant under $\tilde\rho(h)$.  So, there must exist
$g\in\tilde\rho(\pi_1(Y))$ and $x\in\{x_1,x_2\}$ such that
$x_3=g(x)\notin P$, else $\tilde\rho$ would not be irreducible,
again by Lemma \ref{Seifert_1}. As before, we find that $x_3$ is a
fixed point of $\tilde\rho(h)$. By construction, $x_1,x_2$, and
$x_3$ are linearly independent.

Choose a lift of $\tilde\rho(h)$ to ${\rm U}(2,1)$; denote the
lift by $\tilde h$.  The three linearly independent fixed points
$x_1,x_2,x_3$ yield three linearly independent eigenvectors of
$\tilde h$.  We now prove by contradiction that $\tilde h$ has
exactly one eigenvalue.

First, suppose that $\tilde h$ has 3 distinct eigenvalues.  In
this case, we have that $x_1,x_2$, and $x_3$ are exactly the three
one-dimensional eigenspaces of $\tilde h$.  For each
$k\in\{1,\dots,n\}$, lift $\tilde\rho(t_k)$ to ${\rm U}(2,1)$, and
denote the lift by $\tilde t_k$. Now, as before, we find that
$\tilde\rho(t_k)$ maps fixed points of $\tilde\rho(h)$ to fixed
points of $\tilde\rho(h)$.  In other words, $\tilde t_k$ permutes
$x_1,x_2$, and $x_3$.  Let $\eta_k$ be this permutation, regarded
as an element of the symmetric group $S_3$.  The relation
$t_k^{m_k}h^{c_k}=1$ in $\pi_1(Y)$ implies that $\eta_k^{m_k}=1$.
Consequently, the order ${\rm ord}(\eta_k)$ of $\eta_k$ divides
$m_k$.  Now, ${\rm ord}(\eta_k)\in\{1,2,3\}$, and the $m_k$'s are
pairwise coprime.  Therefore, there are at most 2 $k$'s such that
${\rm ord}(\eta_k)\neq 1$.  Moreover, ${\rm ord}(\eta_k)$ is
relatively prime to ${\rm ord}(\eta_{k'})$ whenever $k\neq k'$.
The relation $t_1\dots t_n h^{c_0}=1$ in $\pi_1(Y)$ implies that
$\eta_1\dots\eta_n=1$. Therefore no $\eta_k$ has order 3; for if
so, then $\eta_1\dots\eta_n$ is an odd permutation.  We must then
have that $\eta_k=1$ for each $k$, for otherwise ${\rm
ord}(\eta_1\dots\eta_n)=2$.  However, $\eta_k=1$ if and only if
$\tilde\rho(t_k)$ fixes $x_1,x_2$, and $x_3$.  So every element in
the image of $\tilde\rho$ fixes, say, $x_1$. By Lemma
\ref{Seifert_1}, this contradicts irreducibility of $\tilde\rho$.

Suppose now that $\tilde h$ has exactly 2 distinct eigenvalues.
Without loss of generality, suppose that $x_1$ and $x_2$ belong to
the same 2-dimensional eigenspace $P$ and that $x_3$ is the
1-dimensional eigenspace of $\tilde h$.  Let $f$ be in the image
of $\tilde\rho$, and let $\tilde f$ be a lift of $f$ to ${\rm
U}(2,1)$. We claim that $P$ is invariant under $f$.  As before,
$\tilde f$ commutes with $\tilde h$, so $\tilde f$ maps
eigenvectors of $\tilde h$ to eigenvectors of $\tilde h$.  In
particular, if $P$ is not invariant under $\tilde f$, then $\tilde
f$ maps either $x_1$ or $x_2$ to $x_3$. Let $e_1,e_2$, and $e_3$
be nonzero vectors in $x_1,x_2$, and $x_3$, respectively. Without
loss of generality, assume that $\tilde f(e_1)\in x_3$. Since
$\tilde f$ is nondegenerate, we must then have that $\tilde
f(e_2)\in P$ and $\tilde f(e_3)\in P$. But then $\tilde
f(e_2+e_3)$ is an eigenvector of $\tilde h$ which is neither in
$P$ nor in $x_3$---a contradiction.  So $P$ is invariant under an
arbitrary element in the image of $\tilde\rho$, once again
violating irreducibility.

So, $\tilde h$ has three linearly independent eigenvectors and
exactly one eigenvalue.  Consequently, $\tilde h$ is of the form
$\lambda I$, which implies that $\tilde\rho(h)$ is the identity in
${\rm PU}(2,1)$.$\;\square$

\section{Dolgachev surfaces}
\label{Dolgachev_surfaces}

In this section, we collect facts about our Dolgachev surface $X$
that will be useful later.

Recall the construction of $X$ from \S\ref{Seifert}.  We may
choose our pencil of curves such that each singular fibre is a
rational curve with an ordinary double point. There are, then, 12
such singular fibres in this fibration \cite[p. 192]{Friedman}.
Denote these 12 fibres by $E_1,\dots,E_{12}$. Denote the generic
fibre of $X$ by $F$ and the multiple fibres of $X$ by
$F_{1},\dots,F_{n}$, where $F_{k}$ has multiplicity $m_{k}$.  For
all $j,k$, we have that $E_j$ is linearly equivalent to $F$ is
linearly equivalent to $m_k F_k$.

We say a divisor $D$ on $X$ is \emph{vertical} if $mD$ is linearly
equivalent to $\pi^*(D')$ for some divisor $D'$ on
$\mathbb{CP}^1$.  Note that a multiple fibre $F_k$ is vertical,
but it is not the pullback of a divisor on $\mathbb{CP}^{1}$.
 (Note: this definition of a vertical divisor $D$ is not equivalent
to the condition $D\cdot F=0$, contrary to what one sees
occasionally in the literature.)  A divisor $D$ is vertical if and
only if it is linearly equivalent to $aF + \sum{a_{k}F_{k}}$ for
some integers $a,a_{1},\dots,a_{n}$.  If we write a vertical
divisor in this form, we will always assume that $0\leq a_j<m_j$
for all $j=1,\dots ,n$, unless otherwise noted.

\begin{Lem}[I. Dolgachev]
The surface $X$ has the following numerical invariants: the
topological Euler characteristic $e_X=12$; the irregularity $q=0$;
the geometric genus $p_g=0$.  Also, the canonical bundle
$K_X=\mathcal{O}_X (-F+\sum_k (m_k-1)F_k)$.

\label{Dolgachev_0.5}
\end{Lem}

\begin{Lem}[I. Dolgachev]
$X$ is projective.

\label{Dolgachev_2.5}
\end{Lem}

\begin{Lem}$\;\;\;$\newline
$\;\;\;$(i) $h^0(\mathcal{O}_X(\ell F + \sum\ell_k
F_k))=\rm{max}(\ell +1,0)$, and \newline $\;\;\;$(ii)
$h^1(\mathcal{O}_X(\ell F + \sum\ell_k
F_k))=\rm{max}(\ell,-\ell-1).$

\label{Dolgachev_2}
\end{Lem}
\emph{Proof.} \cite[Lemma 1.1]{BO}

\begin{Lem}
If $s$ is a global section of the locally free sheaf
$\mathcal{O}_X (aF+\sum a_k F_k)$, then $s$ is constant on fibres.

\label{Dolgachev_1}
\end{Lem}
\emph{Proof} Let $w_0=\pi(F)\in\mathbb{CP}^1$.  Choose a local
coordinate $w$ on $\mathbb{CP}^1$ centered at $w_0$.  In order
that $H^0(\mathcal{O}_X (aF+\sum a_k F_k))\neq 0$, we must have
$a\geq 0$, by Lemma \ref{Dolgachev_2}.  Let $f_j=w^{-j}$, for
$j=0,\dots,a$.  The $f_j$'s are linearly independent, so
$\{f_j\circ\pi\}$ is a set of $a+1$ linearly independent elements
in $H^0(\mathcal{O}_X (aF+\sum a_k F_k))\neq 0$.  By Lemma
\ref{Dolgachev_2}, $s$ must be a linearly combination of
$f_j\circ\pi$'s.  Since each $f_j\circ\pi$ is constant on fibres,
so is $s$.$\;\square$

\begin{Lem}
Let $F_k$ be a multiple fibre.  Then there exists a collection
$\{U_{\alpha}\}$ of open sets of $X\!$ such that the
$U_{\alpha}$'s cover $F_k$; each $U_{\alpha}$ is a coordinate
neighborhood on $X$; each $U_{\alpha}$ is disjoint from the
singular fibres and from the other multiple fibres; and, denoting
the coordinates on $U_{\alpha}$ by $(w_{\alpha}, z_{\alpha})$ and
those on $U_{\beta}$ by $(w_{\beta}, z_{\beta})$, we have that
$w_{\alpha}=\zeta_{\alpha\beta}w_{\beta}$ and
$z_{\alpha}=z_{\beta}+t_{\alpha\beta}$ on $U_{\alpha}\cap
U_{\beta}$ for some complex numbers $\zeta_{\alpha\beta}$ with
$\zeta_{\alpha\beta}^{m_k}=1$ and some functions
$t_{\alpha\beta}$; the fibration map $\pi$ locally takes the form
$(w_{\alpha}, z_{\alpha})\stackrel{\pi}{\mapsto}
w=w_{\alpha}^{m_k}$, where $w$ is the local coordinate on
$\mathbb{CP}^{1}$; and $\{w_{\alpha}=0\}$ is a set of local
defining equations for the divisor $F_k$.

\label{Dolgachev_4}
\end{Lem}
\emph{Proof.}  The result follows directly from the definition of
the logarithmic transformation \cite{GH}.  See \cite{me} for more
details.$\;\square$

In the sequel, we will not distinguish between a vector bundle and
its associated locally free sheaf of holomorphic sections, if no
confusion is likely to result.  Two exceptions will come in Lemmas
\ref{Dolgachev_6} and in \S\ref{stable_binary_rank_1}, where we
will make use of the following system of trivializations for
vertical line bundles.

Let $V$ be a small coordinate disc in $\mathbb{CP}^1$, with
coordinate $w$ centered at 0, such that $\pi_0(E_j)\notin V$ for
$j=1,\dots,12$.  Without loss of generality, assume that $V$
contains the points $0, \infty,$ and $\pi(F_k)$ for each multiple
fibre $F_k$; that $\pi(F_k)\notin\{0,\infty\}$ for all $k$; and
that $F=\pi^{-1}(0)$. Cover $\pi^{-1}(V-\infty)-\bigcup F_k$ by
coordinate neighborhoods $V_{\gamma}$ so that there are
coordinates $(w_{\gamma}, z_{\gamma})$ on $V_{\gamma}$, and the
map $\pi$ is given by $\pi(w_{\gamma}, z_{\gamma})=w$ on
$V_{\gamma}$, where $w$ is the coordinate on $\mathbb{CP}^1$
centered at $0$.  For each multiple fibre $F_k$, let $\{U_{\alpha,
k}\}$ be a system of coordinate neighborhoods covering $F_k$,
where $U_{\alpha, k}$ has coordinates $(w_{\alpha, k},z_{\alpha,
k})$.  Cover $\pi^{-1}(V-0)-\bigcup_{\alpha,k}
\overline{U_{\alpha, k}}$ by coordinate neighborhoods $W_{\xi}$ so
that there are coordinates $(w_{\xi}, z_{\xi})$ on $W_{\xi}$, and
the map $\pi$ is given by $\pi(w_{\xi},
z_{\xi})=\frac{1}{w_{\xi}}$ on $W_{\xi}$. The relationships
between the $w$'s are as follows:
$${\rm On}\;U_{\alpha_1}\cap U_{\alpha_2}, {\rm we~have}\;w_{\alpha_1, k}=\zeta_{_{\alpha_1\alpha_2, k}}w_{\alpha_2, k}\;{\rm for~some}\;m_k{\rm th~root~of~unity}\;\zeta_{_{\alpha_1\alpha_2, k}}.$$
$${\rm On}\;U_{\alpha}\cap V_{\gamma}, {\rm we~have}\;w_{\gamma}=w_{\alpha, k}^{m_k}+t_{\alpha, k}\;\;{\rm for~some~complex~number}\;t_{\alpha, k}.$$
$${\rm On}\;V_{\gamma}\cap W_{\xi}, {\rm we~have}\;w_{\xi}=\frac{1}{w_{\gamma}}.$$

Let $L=\mathcal{O}_X (aF+\sum a_k F_k)$ be a vertical line bundle.
Local trivializations for $L$ are given by the maps $f\cdot
w_{\alpha,k}^{-a_k}\mapsto f$ on $U_{\alpha,k}$; $f\cdot
w_{\gamma}^{-a}\mapsto f$ on $V_{\gamma}$; and $f\mapsto f$ on
$W_{\xi}$.  From now on, the notations
$U_{\alpha,k},V_{\gamma},W_{\xi},w_{\alpha,k},w_{\gamma},w_{\xi}$
will be fixed.  Moreover, sections of a vertical line bundle $L$
will be written locally on $U_{\alpha,k},V_{\gamma}$, and
$W_{\xi}$ with respect to these trivializations.

\begin{Lem}
Let $L=\mathcal{O}_X (aF+\sum a_k F_k)$ be a vertical line bundle.

(i) Suppose $a\geq 0$.  If $0\leq j\leq a$, then there exists a
section $s_j\in H^0(L)$ such that $s_j$ is given by
$s_{\xi}=w_{\xi}^{a-j}$ on $W_{\xi}$; $s_{\gamma}=w_{\gamma}^j$ on
$V_{\gamma}$; and $s_{\alpha,k}=(w_{\alpha,k}^{m_k}+t_{\alpha,
k})^j w_{\alpha,k}^{a_k}$ on $U_{\alpha,k}$.  Moreover,
$\{s_j\,|\,0\leq j\leq a\}$ is a basis for $H^0(L)$.

(ii) Suppose $a\leq -2$.  If $a<j<0$, then there exists a \v Cech
1-cocycle $\sigma_j\in C^1(L)$ such that $\sigma_j$ is given by
$\sigma_{\gamma\xi}=w_{\gamma}^j$ on $V_{\gamma}\cap W_{\xi}$ with
respect to the trivialization on $V_{\gamma}$, and
$\sigma_{\xi_1\xi_2}, \sigma_{\gamma_1\gamma_2},
\sigma_{\alpha,k;\gamma}$, and $\sigma_{\alpha_1,k;\alpha_2,k}$
vanish on $W_{\xi_1}\cap W_{\xi_2}, V_{\gamma_1}\cap V_{\gamma_2},
U_{\alpha,k}\cap V_{\gamma}$, and $U_{\alpha_1,k}\cap
U_{\alpha_2,k}$, respectively.  Moreover, identifying $\sigma_j$
with its image in $H^1(L)$, we have that $\{\sigma_j\,|\,a<j<0\}$
is a basis for $H^1(L)$.

\label{Dolgachev_6}
\end{Lem}
\emph{Proof.}  Let $f_j\circ\pi$ be as in the proof of Lemma
\ref{Dolgachev_1}.  Let $s_j=f_j\circ\pi$.  From Lemma
\ref{Dolgachev_1}, we know that $\{s_j\}$ is a basis for $H^0(L)$.
In local coordinates, $s_j$ has the form required in (i).  The
$\sigma_j$'s in (ii) are obtained by pulling back a basis for
$H^1(\mathcal{O}_{\mathbb{CP}^1}(a))$ via $\pi$.$\;\square$

\begin{Def}Let $H_0$ be a fixed ample divisor on $X$.

Let $k_0=1+3\left({\rm
max}\left\{1,-2+\sum\frac{m_k-1}{m_k}\right\}\right)(H_0\cdot F)$.

Let $H=H_0+k_0F$.

\label{Dolgachev_7}
\end{Def}

Note that $H$ is ample. Throughout this document, the degree of a
coherent sheaf---and all related concepts (e.g., stability)---will
be with respect to $H$.

\begin{Lem}
There exists a short exact sequence
$$0 \rightarrow \mathcal{O}_X (-2F+\sum_k (m_k-1)F_k) \rightarrow
\Omega^1_X \rightarrow I_Z\otimes\mathcal{O}_X (F) \rightarrow
0,$$ where $\Omega^1_X$ denotes the sheaf of holomorphic $1$-forms
on $X$, where $Z$ is the reduced subscheme associated to the set
of singular points of singular fibres of $X\!$, and where $I_Z$ is
the ideal sheaf of $Z$.

\label{tern_thm_1}
\end{Lem}
\emph{Proof.}  Pullback of holomorphic 1-forms via $\pi$ gives
rise \cite[p. 98]{BPV} to an injection of sheaves $$0 \rightarrow
\pi^*\Omega^1_{\mathbb{CP}^{1}} \rightarrow \Omega^1_X.$$ Let
$\Omega^1_{X/\mathbb{CP}^{1}}$ denote the sheaf of relative
differentials (i.e., the cokernel of this map). Since
$\pi^*\Omega^1_{\mathbb{CP}^{1}}=\mathcal{O}_{X}(-2F)$, we compute
that ${\rm det}\left( \Omega^1_{X/\mathbb{CP}^{1}}\right) =
\mathcal{O}_{X}(F+\sum (m_k-1)F_k)$.

Let $T={\rm Tor}\left(\Omega^1_{X/\mathbb{CP}^{1}}\right)$, where
${\rm Tor}(\mathcal{S})$ denotes the torsion part of a sheaf
$\mathcal{S}$.  We claim that $T$ is isomorphic to
$\bigoplus_{k=1}^n \mathcal{O}_{(m_k-1)F_k}((m_k-1)F_k)$. To prove
this claim, we first observe that the support of $T$ is contained
in the union of the multiple fibres of $X$ \cite[p. 98]{BPV}.  Let
$F_k$ be a multiple fibre, and let $\{U_{\alpha}\}$ be a
collection of coordinate neighborhoods as in Lemma
\ref{Dolgachev_4}. It suffices to show that $T|_{\cup U_{\alpha}}$
is isomorphic to $\mathcal{O}_{(m_k-1)F_k}((m_k-1)F_k)$.

Let $V$ be an open subset of $\cup U_{\alpha}$.  A section $s$ of
$\Omega^1_{X/\mathbb{CP}^{1}}(V)$ is given by a collection
$\{(V_{\alpha},s_{\alpha})\}$ where $\cup V_{\alpha}=V$,
$s_{\alpha}\in \Omega^1_X(V_{\alpha})$, and
$s_{\beta}-s_{\alpha}\in
\pi^*\Omega^1_{\mathbb{CP}^{1}}(V_{\alpha}\cap V_{\beta})$.
Without loss of generality, we assume that $V_{\alpha}\subset
U_{\alpha}$ for each $\alpha$.  For coordinates on $V_{\alpha}$,
we take the coordinates $(w_{\alpha}, z_{\alpha})$ from
$U_{\alpha}$, as in Lemma \ref{Dolgachev_4}.  Now,
$\Omega^1_X(V_{\alpha})$ is free; its generators are $dw_{\alpha}$
and $dz_{\alpha}$.  Also,
$\pi^*\Omega^1_{\mathbb{CP}^{1}}(V_{\alpha})$ is free, with
generator
$\pi^*(du)=d(w_{\alpha}^{m_k})=(m_k-1)w_{\alpha}^{m_k-1}dw_{\alpha}$,
where $u$ is the local coordinate on $\mathbb{CP}^{1}$. We see
then that locally, $\Omega^1_{X/\mathbb{CP}^{1}}$ has two
generators, $dw_{\alpha}$ and $dz_{\alpha}$, subject to the
relation $w_{\alpha}^{m_k-1}dw_{\alpha}=0$.  Therefore, $T$ is
given locally by the one generator $dw_{\alpha}$ subject to the
relation $w_{\alpha}^{m_k-1}dw_{\alpha}=0$.

Similarly, we find that $\mathcal{O}_{(m_k-1)F_k}((m_k-1)F_k)$ is
given locally by one generator, $w_{\alpha}^{1-m_k}$, subject to
the rather odd-looking relation $w_{\alpha}^{m_k-1}\cdot
w_{\alpha}^{1-m_k}=0$.  Consequently, the map from $T|_{\cup
U_{\alpha}}$ to $\mathcal{O}_{(m_k-1)F_k}((m_k-1)F_k)$ that sends
$dw_{\alpha}$ to $w_{\alpha}^{1-m_k}$ is a well-defined
isomorphism of sheaves.

We can then compute that ${\rm det}(Q)={\rm det}(T)^*\otimes{\rm
det}\left(\Omega^1_{X/\mathbb{CP}^{1}}\right)=\mathcal{O}_{X}(F).$
We have a natural map $\Omega^1_X \rightarrow Q$, which is
surjective. Let $N$ be the kernel of this map.  We then have a
short exact sequence \begin{equation} 0 \rightarrow N \rightarrow
\Omega^1_X \rightarrow Q \rightarrow 0. \label{short_exact_1}
\end{equation}
We then find that $N=\mathcal{O}_X (-2F+\sum(m_k-1)F_k)$. Since
$Q$ is torsion-free, we have that $Q=I_Z\otimes{\rm
det}(Q)=I_Z\otimes\mathcal{O}_X (F)$ for some codimension 2
subscheme $Z$ \cite[p. 33]{Friedman}.  Now,
$\Omega^1_{X/\mathbb{CP}^{1}}$ fails to be locally free precisely
where $\pi$ is singular.  Since $T$ is supported on the union of
the multiple fibres, $Q=\frac{\Omega^1_{X/\mathbb{CP}^{1}}}{T}$
will fail to be locally free at every singular point of $\pi$
outside of the multiple fibres.  In particular, $Z$ contains the
set of singular points of the 12 singular fibres. From
$(\ref{short_exact_1})$ and the equation \cite[p. 29]{Friedman}
$$c_2\left(\Omega^1_X\right)=c_1(N)\cdot c_1\left(\mathcal{O}_X
(F)\right) + \ell(Z)$$ (where $\ell(Z)$ is the length of $Z$), we
find that $\ell(Z)=c_2\left(\Omega^1_X\right)=12$.  We conclude
that $Z$ is the subscheme of $X$ associated to the set of singular
points of the singular fibres, each point taken with multiplicity
one.  The exact sequence (\ref{short_exact_1}) then has the
desired form.$\;\square$

From now on, let $N, Q$, and $Z$ be as in Lemma \ref{tern_thm_1}.

\begin{Lem}
Let $A=aF+\sum a_k F_k$ be a vertical divisor. If
$H^0\!\left(\mathcal{O}_X(-A)\otimes Q\right)\neq 0$, then
$H^0\!\left(\mathcal{O}_X(-A)\otimes N\right)\neq 0$ and
$\rm{deg}(\mathcal{O}_{\it{X}}\!$($A))<0$.

\label{tern_thm_2}
\end{Lem}
\emph{Proof.} A nonzero global section $s$ of
$\mathcal{O}_X(-A)\otimes Q$ is a nonzero global section of
$\mathcal{O}_X(-A+F)$ that vanishes on the total space of $Z$.
Since $-A+F$ is vertical, $s$ is constant on fibres, by Lemma
\ref{Dolgachev_1}. Thus $s$ vanishes identically on each singular
fibre of $X$, and hence can be regarded as a nonzero global
section of $\mathcal{O}_X(-A+F-\sum_{j=1}^{12}(E_j))$.  Now,
$-A+F-\sum_{j=1}^{12}(E_j)$ is linearly equivalent to
$(-11-a-\#\{k|a_k\neq 0\})F+\sum_{a_k\neq 0} (m_k-a_k)F_k$, so by
Lemma \ref{Dolgachev_2}, $a\leq -11-\#\{k|a_k\neq 0\}\leq -2$.
Again by Lemma \ref{Dolgachev_2}, $h^0(\mathcal{O}_X(-A)\otimes
N)=(-2-a)+1>0$, as desired.  Moreover,
$${\rm deg}(\mathcal{O}_X(A))=\left(a+\sum\frac{a_k}{m_k}\right){\rm deg}(F)\leq
\big(a+\#\{k|a_k\neq 0\}\big)\,{\rm deg}(F)\leq -11\, {\rm
deg}(F)<0.\;\square$$

\begin{Lem}
Let $B=bF+\sum b_k F_k$. Then $H^0(\mathcal{O}_X(-B)\otimes
\Omega^1_X)\neq 0$ if and only if $b\leq -2$.

\label{tern_thm_7}
\end{Lem}
\emph{Proof.}  First assume that $b\leq -2$.  Tensoring the exact
sequence (\ref{short_exact_1}) from Lemma \ref{tern_thm_1} with
$\mathcal{O}_X(-B)$, we see that $H^0(\mathcal{O}_X(-B)\otimes
N)\neq 0$. The nonvanishing of $H^0(\mathcal{O}_X(-B)\otimes N)$
follows from the effectiveness of $-B+(-2F+\sum_k (m_k-1)F_k).$
(Recall the convention that $b_k<m_k$ for all $k$.)

We now assume that $H^0(\mathcal{O}_X(-B)\otimes \Omega^1_X)\neq
0$ and show that $b\leq -2$.  We must have
$H^0(\mathcal{O}_X(-B)\,\otimes\, Q)\neq 0$ or
$H^0(\mathcal{O}_X(-B)\,\otimes\, N)\neq 0$. Either way,
$H^0(\mathcal{O}_X(-B)\,\otimes\, N)\neq 0$, by
Lemma~\ref{tern_thm_2}.  But then $(-2-b)F+\sum (m_k-1-b_k)F_k$ is
linearly equivalent to an effective divisor.  Therefore $b\leq
-2$.$\;\square$

\textbf{Remark:}  In fact, we can compute that
$h^0(\mathcal{O}_X(-B)\otimes \Omega^1_X)={\rm max}\{0,-2-b\}$. To
do so, let $L=\mathcal{O}_X(B)$, and consider the exact sequence
$0\to\pi_*(L^{^*}\!\otimes
N)\to\pi_*(L^{^*}\!\otimes\Omega^1_X)\to\pi_*(L^{^*}\!\otimes Q).$
Then show that $\pi_*(L^{^*}\!\otimes N)$ is a line bundle on
$\mathbb{CP}^1$, that $\pi_*(L^{^*}\!\otimes\Omega^1_X)$ is a
coherent sheaf of rank 1 on $\mathbb{CP}^1$, and that
$\pi_*(L^{^*}\!\otimes Q)$ is torsion-free.  It follows that
$${\rm max}\{0,-2-b\}=h^0(L^{^*}\!\otimes N)=h^0(\pi_*(L^{^*}\!\otimes
N))=h^0(\pi_*(L^{^*}\!\otimes\Omega^1_X))=h^0(L^{^*}\!\otimes\Omega^1_X).$$

\section{Toledo invariants on 2-orbifolds and Dolgachev surfaces}
\label{orbifold}

In this section, we associate to our Seifert fibered space $Y$ a
2-orbifold $O$.  The goal of this section is to show how Toledo
invariants on the Dolgachev surface $X$ correspond to ``orbifold''
Toledo invariants which arise from representations of the orbifold
fundamental group of $O$.

Let $O$ be the hyperbolic 2-orbifold such that the underlying
space $|O|$ of $O$ is the sphere $S^2$ and $O$ has $n$ elliptic
points $p_1,\dots,p_n$ (also known as cone points) of orders
$m_1,\dots,m_n$, respectively. (We refer to \cite{Boden},
\cite{FS}, \cite{Kapovich}, \cite{Scott}, and \cite{Thurston} for
details of this construction and for basic facts about orbifolds.)
The orbifold fundamental group of $O$ has the following
presentation:
$$\pi_1^{\rm{orb}}(O)=<u_1,\dots,u_n\,|\,u_k^{m_k}=u_1\dots u_n=1>$$

We may think of $u_j$ as a small loop that travels once around the
cone point $p_j$.

In our elliptic fibration $\pi:X\rightarrow\mathbb{CP}^1$, we
identify $\mathbb{CP}^1$ with $|O|$, and we assume that
$p_j=\pi(F_j)$ for each multiple fibre $F_j$.  Let $\tilde X$ be
the universal cover of our Dolgachev surface $X$.  The
restrictions we imposed on the $m$'s following Lemma
\ref{Seifert_0.5} imply that the orbifold universal cover $\tilde
O$ of $O$ is the upper half-plane $H^2$ \cite{Thurston}. Fix a
base point $x_0$ in $X$ and a base point $y_0$ in $O$ such that
$y_0=\pi(x_0)$ and $x_0\notin\{ E_1,\dots,E_{12},F_1,\dots, F_n
\}$. We may regard the elements of $\tilde X$ (resp. $\tilde O$)
as equivalence classes of paths in $X$ (resp. $O$) beginning at
$x_0$ (resp. $y_0$). (Caution: One must be careful as to what is
meant by a path in $O$.  See \cite[\S 2]{FS}.) Pushing forward
paths in $X$ to paths in $O$, we obtain a map $\tilde\pi:\tilde
X\rightarrow\tilde O$ that covers $\pi$.  If $\gamma$ is an
element of $\pi_1(X)$, then denote the action of $\gamma$ on
$\tilde X$ by $L_{\gamma}$. (Similarly for $\tilde O$.) Recall
that $t_1,\dots,t_n$ are the generators of $\pi_1(X)$. Then
$\tilde\pi\circ L_{t_j}=L_{u_j}\circ\tilde\pi$. It follows that
$\pi_*(t_j)=u_j$, and so $\pi_*$ is an isomorphism from $\pi_1(X)$
to $\pi_1^{\rm{orb}}(O)$.

\begin{Lem}
Let $\rho\in{\rm Hom}(\pi_1(X),{\rm U}(2,1))$, and let
$E_{\rho}=\frac{\tilde X\times\emph{\scriptsize\textbf
H}_{\mathbb{C}}^2}{\pi_1(X)}$. Let $q:\tilde
X\times\emph{\textbf{H}}_{\mathbb{C}}^2\rightarrow\emph{\textbf{H}}_{\mathbb{C}}^2$
be projection onto the second factor.  Then there exists a section
$s_0$ of the fibre bundle $E_{\rho}$ (as in Definition {\rm
\ref{Toledo_0.5}}) and a lift $\tilde s_0:\tilde
X\rightarrow\tilde X\times\emph{\textbf{H}}_{\mathbb{C}}^2$ of
$s_0$ such that for each point $x\in\tilde O$, we have that
$q\circ\tilde s$ is constant on $\tilde\pi^{-1}(x)$.

\label{orbifold_0.5}
\end{Lem}

\begin{Def}
Let $\rho\in{\rm Hom}(\pi_1^{\rm{orb}}(O),\rm{PU}(2,1))$.  We say
that a map $s:\tilde O\rightarrow\tilde
O\times\emph{\textbf{H}}_{\mathbb{C}}^2$ is
\emph{$\rho$-equivariant} if $s(\gamma\cdot x)=\rho(\gamma)\cdot
s(x)$ for all $\gamma\in\pi_1^{\rm{orb}}(O)$ and $x\in\tilde O$.
If $s_1$ and $s_2$ are two $\rho$-equivariant maps, then we say
that $s_1$ and $s_2$ are \emph{$\rho$-equivariantly homotopic} if
there exists a homotopy $F:[0,1]\times\tilde
O\rightarrow[0,1]\times\tilde
O\times\emph{\textbf{H}}_{\mathbb{C}}^2$ from $s_1$ to $s_2$ such
that $F(t,\cdot,\cdot)$ is $\rho$-equivariant for all $t\in[0,1]$.

\label{orbifold_1}
\end{Def}

\begin{Lem}
Let $\rho\in{\rm Hom}(\pi_1^{\rm{orb}}(O),\rm{PU}(2,1))$.  Then
there exists a $\rho$-equivariant map\linebreak $s:\tilde
O\rightarrow\tilde O\times\emph{\textbf{H}}_{\mathbb{C}}^2$.
Moreover, if $s_1$ and $s_2$ are any two such $\rho$-equivariant
maps, then $s_1$ and $s_2$ are $\rho$-equivariantly homotopic.

\label{orbifold_2}
\end{Lem}
\emph{Proof.}  Existence: push forward $\tilde s_0$ from Lemma
\ref{orbifold_0.5}.  Invariance: push forward a homotopy.

\begin{Def}
Let $\rho\in{\rm Hom}(\pi_1^{\rm{orb}}(O),\rm{PU}(2,1))$.  Let
$\Sigma$ be a fundamental domain for the action of
$\pi_1^{\rm{orb}}(O)$ on $\tilde O$. Take $q$ and $\omega$ as in
Lemma {\rm \ref{orbifold_0.5}}. Then we define the \emph{orbifold
Toledo invariant} $\tau_{\rm{orb}}(\rho)$ by
$$\tau_{\rm{orb}}(\rho)=\int_{\Sigma}s^*q^*\omega$$where $s:\tilde O\rightarrow\tilde
O\times\emph{\textbf{H}}_{\mathbb{C}}^2$ is any $\rho$-equivariant
map.

\label{orbifold_3}
\end{Def}

Lemma \ref{orbifold_2} implies that $\tau_{\rm{orb}}(\rho)$ is
defined and that it is independent of the choice of $s$.  The
$\rho$-equivariance of $s$ implies that $\tau_{\rm{orb}}(\rho)$ is
independent of the choice of $\Sigma$.  We now fix $\tilde s_0$ as
in Lemma \ref{orbifold_0.5}, and let $s$ be its $\rho$-equivariant
push-forward, as in Lemma \ref{orbifold_2}.

Let $H^2_{\rm{orb}}(O,\mathbb{Z})$ be the orbifold second
cohomology group of $O$ with integer coefficients \cite{FS}. (Note
that \cite{FS} uses the notation ``$V$'' in place of ``orb,''
since they use the older terminology ``V-manifold'' in place of
``orbifold.'')  Let $H^1_{\rm{vert}}(X,\mathcal{O}_X^{^*}\!)$ be
the subgroup of $H^1(X,\mathcal{O}_X^{^*}\!)$ consisting of
vertical line bundles on $X$. Let
$H^2_{\rm{vert}}(X,\mathbb{Z})=c_1(H^1_{\rm{vert}}(X,\mathcal{O}_X^{^*}\!))$
be the group of first Chern classes of vertical line bundles on
$X$.

\begin{Lem}
The map $\pi$ induces an isomorphism
$\pi^*:H^2_{\rm{orb}}(O,\mathbb{Z})\rightarrow
H^2_{\rm{vert}}(X,\mathbb{Z})$.

\label{orbifold_5}
\end{Lem}

Let ${\rm Pic}^t_{\rm{orb}}(O)$ be the set of topological
isomorphism classes of orbifold line bundles on $O$.  We have that
${\rm Pic}^t_{\rm{orb}}(O)$ is a group, where the group law is
given by the tensor product.

\begin{Lem}
${\rm Pic}^t_{\rm{orb}}(O)\cong H^2_{\rm{orb}}(O,\mathbb{Z})$

\label{orbifold_6}
\end{Lem}
\emph{Proof.} \cite[Theorem 2.2(ii)]{FS}$\;\;\;\square$

\begin{Lem}
Let $\rho\in{\rm Hom}(\pi_1^{\rm{orb}}(O),\rm{PU}(2,1))$.  If
$\tau(\rho\circ\pi_*)=c_1(\mathcal{O}_X(aF+\sum a_kF_k))$, then
$\tau_{\rm{orb}}(\rho)=a+\sum\frac{a_k}{m_k}$.

\label{orbifold_7}
\end{Lem}
\emph{Proof.} Let $p\in|O|-\{p_1,\dots,p_n\}$. Let $L_p$ be the
holomorphic point bundle determined by $p$.  Then $L_p$ is an
orbifold line bundle on $O$ with $c_1(L_p)=1$ \cite{FS}.  Let
$\sigma_k:\frac{\mathbb{Z}}{m_k\mathbb{Z}}\rightarrow U(1)$ be the
standard representation.  Let $L_{p_k}$ be the orbifold line
bundle on $O$ with first Chern class $c_1(L_{p_k})=\frac{1}{m_k}$
and trivial isotropy except at $p_k$, where it is $\sigma_k$. Then
$\pi^*L_p=\mathcal{O}_X(F)$ and $\pi^*L_{p_k}=\mathcal{O}_X(F_k)$.
Let $L=L_p^{\otimes a}\otimes\left(\bigotimes_k L_{p_k}^{\otimes
a_k}\right)$.  The following diagrams commute:

\begin{tabular}{ccccccc}

& & $\;\;\;\tilde X$ & $\stackrel{\tilde s_0}{\longrightarrow}$ &
$\tilde X\times\textbf{H}_{\mathbb{C}}^2$ & $\stackrel{q_X}{\longrightarrow}$ & $\textbf{H}_{\mathbb{C}}^2$ \vspace*{.05in} \\

& & $\tilde\pi\downarrow$        &                                               & $\downarrow$ \vspace*{.05in} & & $\downarrow$\\

& & $\;\;\;\tilde O$ & $\stackrel{s}{\longrightarrow}$ &
$\tilde O\times\textbf{H}_{\mathbb{C}}^2$ & $\stackrel{q_O}{\longrightarrow}$ & $\textbf{H}_{\mathbb{C}}^2$ \vspace*{.05in} \\

\end{tabular}

and

\begin{tabular}{ccccc}

& & $\;\;\;\tilde X$ & $\stackrel{\tilde\pi}{\longrightarrow}$ &
$\tilde O$ \vspace*{.05in} \\

& & $\varphi_X\downarrow$        &                                               & $\downarrow\varphi_O$ \vspace*{.05in} \\

& & $\;\;\;X$ & $\stackrel{\pi}{\longrightarrow}$ &
$O$ \vspace*{.05in} \\

\end{tabular}

and

\begin{tabular}{ccccc}

& & $\;\;\;\tilde X$ & $\stackrel{\tilde s_0}\longrightarrow$ &
$\tilde X\times\textbf{H}_{\mathbb{C}}^2$ \vspace*{.05in} \\

& & $\varphi_X\downarrow$        &                                               & $\downarrow\varphi$ \vspace*{.05in} \\

& & $\;\;\;X$ & $\stackrel{s_0}\longrightarrow$ &
$E_{\rho\circ\pi_*}$ \vspace*{.05in} \\

\end{tabular}

From these diagrams, we find that $\tilde\pi^*c_1(\varphi_O^*L)
=c_1(\tilde\pi^*\varphi_O^*L) =c_1(\varphi_X^*\pi^*L)=\linebreak
\varphi_X^*c_1(\mathcal{O}_X(aF+\sum a_kF_k)
=\varphi_X^*\tau(\rho\circ\pi_*)=\varphi_X^*s_0^*(\varphi_X)_*q_X^*\omega=\tilde
s_0^*q_X^*\omega=\tilde\pi^*s^*q_O^*\omega$.

Let us identify $H^2_{\rm{orb}}(O,\mathbb{Z})$ with the set of all
$\pi_1^{\rm{orb}}(O)$-invariant elements of $H^2(\tilde
O,\mathbb{Z})$. Lemma \ref{orbifold_5} then implies that
$c_1(\varphi_O^*L)=s^*q_O^*\omega$.  Therefore
$$\tau_{\rm{orb}}(\rho)=\int_{\Sigma}s^*q_O^*\omega=\int_{\Sigma}c_1(\varphi_O^*L)=c_1(L)=a+\sum\frac{a_k}{m_k}.\;\;\;\;\square$$

\section{U(2,1) Higgs bundles}
\label{Higgs}

Hitchin, Simpson, et al. (\cite{Hitchin}, \cite{Simpson},
\cite{Simpson2}) have shown that representations of the
fundamental group of a compact K\"ahler manifold are closely
related to holomorphic objects called Higgs bundles.  The goal of
this section is to describe the Higgs bundles that arise from
U(2,1) representations of the fundamental group of a Dolgachev
surface, then describe the Toledo invariant of such a
representation in terms of the Chern classes of the associated
Higgs bundle.

\begin{Def}
Let $M$ be a complex algebraic manifold, and let $H$ be a fixed
ample line bundle on $M$. A \emph{Higgs bundle} on $M$ is a pair
$(V,\theta)$, where $V$ is a holomorphic vector bundle on $M$;
$\theta\!\in H^0($\emph{End}$(V)\otimes\Omega^1_M)$; and
$\theta\wedge\theta=0$.  $\theta$ is called the \emph{Higgs
field}. A subsheaf $\mathcal{S}$ of $V$ is said to be
\emph{$\theta$-invariant} if
$\theta(\mathcal{S})\subset\mathcal{S}\otimes\Omega^1_M$. The
\emph{slope} $\mu(\mathcal{S})$ of a coherent sheaf $\mathcal{S}$
on $M$ with \emph{rank}$(\mathcal{S})>0$ is defined by
$\mu(\mathcal{S})=\frac{\rm{deg}(\mathcal{S})}{\rm{rank}(\mathcal{S})}$,
where \emph{deg}$(\mathcal{S})$ is the degree of $\mathcal{S}$
with respect to $H$.  A Higgs bundle $(V,\theta)$ is \emph{stable}
if $\mu(\mathcal{S})<\mu(V)$ for all coherent $\theta$-invariant
subsheaves $\mathcal{S}$ of $V$ with \emph{rank}$(\mathcal{S})>0$.
 A Higgs bundle $(V,\theta)$ is \emph{polystable} if it is a direct
sum of stable Higgs bundles, each with the same slope.  (One forms
the direct sum in the obvious way.)  A Higgs bundle $(V,\theta)$
is \emph{reducible} if it is a direct sum of Higgs bundles and is
\emph{irreducible} otherwise.  We say that a Higgs bundle
$(V,\theta)$ is a \emph{U(2,1)-Higgs bundle} if $V = V_P\oplus
V_Q$ (where $V_P$ and $V_Q$ are vector bundles of rank 2 and 1,
respectively), and $\theta$ maps $V_P$ to $V_Q\otimes\Omega^1_M$
and $V_Q$ to $V_P\otimes\Omega^1_M$.

\label{Higgs_1}
\end{Def}

If $H$ is any group, then let ${\rm Hom}^+(H,\rm{U}(2,1))$ denote
the space of semisimple representations from $H$ into
$\rm{U}(2,1)$.

\begin{Lem}
There exists a surjective function
$\phi:\mathcal{H}\rightarrow\;$\emph{Hom}$^+(\pi_1(X),\rm{U}(2,1))$,
where $\mathcal{H}$ is the set of all polystable \emph{U(2,1)}
Higgs bundles $(V,\theta)$ on $X$ whose summands have vanishing
Chern classes.

\label{Higgs_2}
\end{Lem}
\emph{Proof.}  Let $\mathcal{H}'$ be the set of all polystable
rank 3 Higgs bundles $(V,\theta)$ on $X$ whose summands have
vanishing Chern classes. By Lemma \ref{Dolgachev_2.5}, $X$ is
algebraic, hence compact K\"ahler. In \cite{Simpson}, Simpson
shows that there is a surjective function
$\phi:\mathcal{H}'\rightarrow{\rm
Hom}^+(\pi_1(X),\rm{GL}(3,\mathbb{C}))$. In \cite[Proposition
3.1]{Xia}, Xia shows that $\mathcal{H}=\phi^{-1}({\rm
Hom}^+(\pi_1(X),\rm{U}(2,1)))$. (Xia's proof is for Riemann
surfaces, but it goes through for any compact K\"ahler
manifold.)$\;\square$

\begin{Lem}
Let $\mathcal{H}$ and $\phi$ be as in Lemma {\rm \ref{Higgs_2}},
and let $(V,\theta)\in\mathcal{H}$.  Write $V=V_P\oplus V_Q$ as in
Def. {\rm \ref{Higgs_1}}.  Then $\tau(\phi(V,\theta))=c_1(V_P)$.

\label{Higgs_3}
\end{Lem}
\emph{Proof.}  \cite{Xia}

\textbf{Remark:} Lemmas \ref{Higgs_2} and \ref{Higgs_3} serve as a
bridge from the world of semisimple representations and Toledo
invariants to the world of polystable Higgs bundles and Chern
classes.  Consequently, while the definition of the Toledo
invariant is topological in nature, these two lemmas enable us to
use algebraic geometry in order to compute which Toledo invariants
actually occur.

\begin{Def}
If $(V,\theta)=(V_P\oplus V_Q,\theta)$ is a $\rm{U}(2,1)$ Higgs
bundle as in Def. {\rm \ref{Higgs_1}}, then we define the
\emph{Higgs bundle Toledo invariant} $\tau_{(V,\theta)}$ by
$\tau_{(V,\theta)}=\frac{1}{3}(c_1(V_P)-2c_1(V_Q))$.

\label{Higgs_4}
\end{Def}

Note that if $(V,\theta)\in\mathcal{H}$, then $V$ is flat, in
which case Def. \ref{Higgs_4} is consistent with Lemma
\ref{Higgs_3}.

\textbf{Remark:}  Consider a semisimple representation
$\rho:\pi_1(X)\to\rm{PU}(2,1)$.  The corresponding principal
$\rm{PU}(2,1)$ bundle on $X$ lifts to a principal $\rm{U}(2,1)$
bundle with an associated vector bundle $V=V_P\oplus V_Q$.  One
can show that $\tau(\rho)=\frac{1}{3}(c_1(V_P)-2c_1(V_Q))$; the
proof is similar to that of Lemma \ref{Higgs_3}.  This is the
motivation for Def. \ref{Higgs_4}.

\begin{Lem}
Let $(V,\theta)\in\mathcal{H}$, and let $L$ be a line bundle.
Then:

{\rm (i)} $(V\otimes L,\theta\otimes 1)$ is a polystable
$\rm{U}(2,1)$ Higgs bundle with $\tau_{(V\otimes L,\,\theta\otimes
1)}=\tau_{(V,\,\theta)}$.

{\rm (ii)} $(V^{^*}\!,\theta)\in\mathcal{H}$, and
$\tau_{(V^{^*}\!,\,\theta)}=-\tau_{(V,\,\theta)}$.

\label{Higgs_5}
\end{Lem}
\emph{Proof.}  These statements follow directly from the
definitions.  See \cite{Xia} for more details.

\section{Systems of Hodge bundles on Dolgachev surfaces}
\label{Hodge}

The results of \S\ref{Higgs} imply that to compute Toledo
invariants of semisimple ${\rm U}(2,1)$ representations of the
fundamental group of a Dolgachev surface, it suffices to compute
Chern classes of the summands of certain polystable ${\rm U}(2,1)$
Higgs bundles.  The goal of this section is to show that we may
restrict our attention to a special class of these Higgs bundles,
namely systems of Hodge bundles.  The method is due to Simpson.
(See \cite{Simpson} and \cite{Simpson3}).  Following Xia
\cite{Xia}, we then divide these systems of Hodge bundles into two
types, binary and ternary.

\begin{Def}[\cite{Simpson3}]
Let $M$ be a complex algebraic manifold.  A \emph{system of Hodge
bundles} on $M$ is a Higgs bundle $(V,\theta)$ such that
$V=\bigoplus V^{r,s}$ and $\theta:V^{r,s}\rightarrow
V^{r-1,s+1}\otimes\Omega^1_M$.

\label{Hodge_1}
\end{Def}

\begin{Lem}
{\rm (a)} There exists a quasiprojective variety
$\mathcal{M}_{\rm{Dol}}$ whose points parametrize direct sums of
stable Higgs bundles with vanishing Chern classes.

{\rm (b)} Let $f$ be the map from $\mathcal{M}_{\rm{Dol}}$ to the
space of polynomials with coefficients in symmetric powers of the
cotangent bundle which takes $(V,\theta)$ to the characteristic
polynomial of $\theta$.  Then $f$ is proper.

{\rm (c)} Let $\mathcal{M}_{\rm{Dol}}(\rm{U}(2,1))$ denote the
subspace of $\mathcal{M}_{\rm{Dol}}$ whose points parametrize
polystable \emph{U(2,1)} bundles.  Then every connected component
of $\mathcal{M}_{\rm{Dol}}(\rm{U}(2,1))$ contains a system of
Hodge bundles.

{\rm (d)} $\mathcal{M}_{\rm{Dol}}(\rm{U}(2,1))$ is homeomorphic to
$\mathcal{R}^+_{\rm{U}(2,1)}(X)$.

\label{Hodge_2}
\end{Lem}
\emph{Proof.}$\;$(a) \cite[Prop. 1.4]{Simpson3}

(b) \cite[Prop. 1.4]{Simpson3}

(c) In \cite[Theorem 3]{Simpson3}, Simpson proves that every
component of $\mathcal{M}_{\rm{Dol}}$ contains a system of Hodge
bundles, as follows.  Let $\mathbb{C}^*$ act on
$\mathcal{M}_{\rm{Dol}}$ by $t\cdot(V,\theta)=(V,t\theta)$.  As
$t\rightarrow 0$, we have $f(t\cdot(V,\theta))\rightarrow 0$.
Since $f$ is proper, $t\cdot(V,\theta)$ converges to a limit Higgs
bundle $(V_0,\theta_0)$.  Since $\mathcal{M}_{\rm{Dol}}$ is
Hausdorff, the limit is unique.  Consequently, $(V_0,\theta_0)$ is
fixed under the action of $\mathbb{C}^*$ and is therefore a system
of Hodge bundles \cite[Lemma 4.1]{Simpson3}.

Since ${\rm U}(2,1)$ is closed in GL$(3,\mathbb{C})$, we have that
$\mathcal{M}_{\rm{Dol}}(\rm{U}(2,1))$ is closed in
$\mathcal{M}_{\rm{Dol}}$.  Therefore, $f$ restricted to
$\mathcal{M}_{\rm{Dol}}(\rm{U}(2,1))$ is still proper, and the
above proof goes through unchanged.$\;\square$

\begin{Def}[\cite{Xia}]
We say a Higgs bundle $(V,\theta)$ is \emph{binary} if $V =
V_P\oplus V_Q$ where $V_P$ and $V_Q$ are vector bundles of rank
$2$ and $1$, respectively; and $\theta$ maps $V_P$ to
$V_Q\otimes\Omega^1_X$ and $V_Q$ to $0$.  In this situation,
denote $(V,\theta)$ by $V_P
\stackrel{\theta\oplus}{\longrightarrow} V_Q$ (omitting $\theta$
if it's clear from the context).

We say a Higgs bundle $(V,\theta)$ is \emph{ternary} if $V =
V_2\oplus V_3\oplus V_1$ where $V_1$, $V_2$, and $V_3$ are line
bundles; and $\theta$ maps $V_2$ to $V_3\otimes\Omega^1_X$, maps
$V_3$ to $V_1\otimes\Omega^1_X$, and maps $V_1$ to $0$.  In this
situation, denote $(V,\theta)$ by
$V_2\stackrel{\oplus}{\rightarrow}
V_3\stackrel{\oplus}{\rightarrow} V_1$.  (In this case, we have
$V_P=V_1\oplus V_2$ and $V_Q=V_3$.)

\label{Hodge_3}
\end{Def}
It follows from Def. \!\ref{Hodge_3} and Lemma \ref{Higgs_5} that
if a polystable ${\rm U}(2,1)$ Higgs bundle is a system of Hodge
bundles, then it is either ternary, binary, or dual to a binary
bundle.  Also, every polystable Higgs bundle is either stable or
reducible.  We therefore investigate the following four types of
polystable ${\rm U}(2,1)$ Higgs bundles: stable ternary, stable
binary, reducible ternary, and reducible binary.

\subsection{The case of the stable ternary Higgs bundle}

\begin{Prop}
Let $V_2=\mathcal{O}_X(bF+\sum{b_{k}F_{k}})$ and
$V_1=\mathcal{O}_X(aF+\sum{a_{k}F_{k}})$.  Then there exists a
Higgs field $\theta$ such that
$(V,\theta)=V_2\stackrel{\oplus}{\rightarrow}\mathcal{O}_X\stackrel{\oplus}{\rightarrow}
V_1$ is a stable ternary Higgs bundle if and only if:

\emph{(i)} $b\leq -2$, and

\emph{(ii)} $a+\#\{k\thinspace |\thinspace a_k\neq 0\}\geq 2$, and

\emph{(iii)} $2A<B$, and

\emph{(iv)} $A<2B$.

Here we have used the notations $A=a+\sum{\frac{a_k}{m_k}}$ and
$B=b+\sum{\frac{b_k}{m_k}}$.

\label{tern_thm_11}
\end{Prop}

\emph{Proof}.  First assume that such a Higgs field $\theta$
exists.  Stability then implies that $\theta|_{V_2}$ and
$\theta|_{\mathcal{O}_X}$ are nonzero. Hence
$H^0(V_2^{^*}\!\otimes\Omega^1_X)\neq 0$ and
$H^0(V_1\otimes\Omega^1_X)\neq 0$.  Conditions (i) and (ii) then
follow from Lemma \ref{tern_thm_7}.  Conditions (iii) and (iv)
follow from the fact that the $\theta$-invariant subsheaves
$\mathcal{O}_X\oplus V_1$ and $V_1$ are not destabilizing.

Conversely, if (i) and (ii) hold, then let $\theta_2$ be a nonzero
global map from $V_2$ to $N$ and $\theta_1$ a nonzero global map
from $\mathcal{O}_X$ to $V_1\otimes N$.  (Lemma \ref{tern_thm_7}
shows that $\theta_2$ and $\theta_1$ exist.) Let
$(w_{\gamma},z_{\gamma})$ be coordinates on $V_{\gamma}$, as in
the discussion following Lemma \ref{Dolgachev_4}.  On
$V_{\gamma}$, then, $\theta_1$ has the form $g_1 dw_{\gamma}$ for
some meromorphic function $g_1$, and $\theta_2=g_2 dw_{\gamma}$ on
$V_{\gamma}$ for some meromorphic $g_2$.  Define $\theta$ by
$\theta|V_2=\theta_2$, $\theta|\mathcal{O}_X=\theta_1$, and
$\theta|V_1=0$.  Then
$\theta\wedge\theta=\theta_1\wedge\theta_2=0$ on $V_{\gamma}$.
Similarly, we find that $\theta\wedge\theta$ vanishes outside the
union of the singular fibres and the multiple fibres.  Hence
$\theta\wedge\theta=0$ everywhere.  Moreover, conditions (iii) and
(iv), together with the nonvanishing of $\theta_1$ and $\theta_2$,
guarantee that $(V,\theta)$ is stable.$\;\square$

\begin{Lem}
Suppose that
$(V,\theta)=V_2\stackrel{\oplus}{\rightarrow}\mathcal{O}_X\stackrel{\oplus}{\rightarrow}
V_1$ is a stable ternary Higgs bundle.  Then $V_2$ and $V_1$ are
vertical.

\label{tern_thm_9}
\end{Lem}
\emph{Proof.}  Choose divisors $D_1$ and $D_2$ such that
$V_1=\mathcal{O}_X(D_1)$ and $V_2=\mathcal{O}_X(D_2)$.  As in the
proof of Lemma \ref{tern_thm_11}, we see that
$H^0(\mathcal{O}_X(-D_2)\otimes\Omega^1_X)\neq 0$.  From the short
exact sequence (\ref{short_exact_1}) in Lemma \ref{tern_thm_1}, we
find that either $-D_2-2F+\sum(m_k-1)F_k$ or $-D_2+F$ is linearly
equivalent to an effective divisor.  Consequently, we find that
$D_2\cdot F\leq 0$, with equality iff $D_2$ is vertical.  We also
find that $H_0\cdot D_2\leq\frac{k_0}{3}$, where $H_0$ and $k_0$
are as in Def. \ref{Dolgachev_7}.  Similarly, we find that
$D_1\cdot F\geq 0$, with equality iff $D_1$ is vertical, and that
$H_0\cdot D_1\geq -\frac{k_0}{3}$.  Therefore
$H_0\cdot(D_1-2D_2)\geq -k_0$. Suppose that either $D_1$ or $D_2$
is nonvertical.  Then, from Def. \ref{Dolgachev_7}, we find that
$H\cdot(D_1-2D_2)\geq 0$.  But this violates (iv) of Prop.
\ref{tern_thm_11}.$\;\square$

\subsection{The case of the stable binary Higgs bundle with rank(im($\theta$))=1}
\label{stable_binary_rank_1}

~\newline Let
$(V,\theta)=V_P\stackrel{\oplus}{\longrightarrow}\mathcal{O}_X$ be
a stable projectively flat binary Higgs bundle. When restricted to
$V_P$, the Higgs field $\theta|V_P$ is a map from $V_P$ to
$\Omega^1_X$.  The image ${\rm im}(\theta|V_P)$ of this map is a
subsheaf of $\Omega^1_X$.  Stability implies that $\theta|V_P$
cannot be the zero map.  It follows that ${\rm im}(\theta|V_P)$
has rank 1 or rank 2.  We shall take these cases separately,
beginning with the rank 1 case.

\begin{Prop}
If $(V,\theta)=V_P\stackrel{\oplus}{\longrightarrow}\mathcal{O}_X$
is a stable projectively flat binary Higgs bundle with
$\emph{rank}(\emph{im}(\theta|V_P))=1$, then $V_P$ can be written
as an extension of the form \begin{equation}0\rightarrow
V_1\rightarrow V_P\stackrel{\beta}{\rightarrow} V_2 \rightarrow
0,\label{short_exact_4}\end{equation} where
$V_1=\mathcal{O}_X(aF+\sum a_k F_k)$ and
$V_2=\mathcal{O}_X(bF+\sum b_k F_k)$ with the $a$'s and $b$'s
subject to the following numerical conditions:

\begin{itemize}

\item[(i)] $-B<A<\frac{1}{2}B$, and \\

\vspace{-.2in}

\item[(ii)] $d_2\leq -2$, and \\

\vspace{-.2in}

\item[(iii)] $b\leq -2$, and \\

\vspace{-.2in}

\item[(iv)] If $(c,c_1,\dots,c_n)$ is an $(n+1)$-tuple of integers
such that $0\leq c_k<m_k$ for all $k$ and $d_1\geq 0$ and
$C\geq\frac{2}{3}(A+B)$, then
$d_1+1\leq\emph{min}(-d_2-1,-d_3-1)$.

\end{itemize}

Here we have used the notations $A=a+\sum\frac{a_k}{m_k}$;
$B=b+\sum\frac{b_k}{m_k}$; $C=c+\sum\frac{c_k}{m_k}$;
$d_1=b-c-\#\{b_k<c_k\}$; $d_2=a-b-\#\{a_k<b_k\}$; and
$d_3=a-c-\#\{a_k<c_k\}$.

Conversely, given $a$'s and $b$'s satisfying \emph{(i)--(iv)},
there exists a stable projectively flat binary Higgs bundle
$V_P\stackrel{\oplus}{\longrightarrow}\mathcal{O}_X$ with $V_P$
given as an extension of the form \rm{(\ref{short_exact_4})}.

\label{bin_thm_20}
\end{Prop}
Before proving this proposition, we first prove several
preliminary lemmas.

\begin{Lem}
Let
$(V,\theta)=V_P\stackrel{\oplus}{\longrightarrow}\mathcal{O}_X$ be
a binary Higgs bundle such that ${\rm im}(\theta|V_P)$ has rank 1.
Let $V_1={\rm ker}(\theta|V_P)$.  Then $(V,\theta)$ is stable if
and only if:

{\rm (SB1)} ${\rm deg}(V_1)<\frac{1}{3}\,{\rm deg}(V_P)$, and

{\rm (SB2)} ${\rm deg}(\mathcal{S})<\frac{2}{3}\,{\rm deg}(V_P)$
for every rank 1 subsheaf $\mathcal{S}$ of $V_P$, and

{\rm (SB3)} ${\rm deg}(V_P)>0.$

\label{bin_thm_1}
\end{Lem}
\emph{Proof.}  If $(V,\theta)$ is stable, then (SB1)--(SB3) follow
directly from the fact that the $\theta$\verb'-'invariant
subsheaves $V_1, \mathcal{S}\oplus\mathcal{O}_X$, and
$\mathcal{O}_X$, respectively, do not destabilize $V$. Conversely,
if \hbox{(SB1)--(SB3)} hold, then any proper $\theta$-invariant
subsheaf $\mathcal{S}'$ of $V$ must be a \hbox{rank 1} subsheaf of
$V_1$, a \hbox{rank 1} subsheaf of $\mathcal{O}_X$, or of the form
$\mathcal{S}\oplus\mathcal{O}_X$, where $\mathcal{S}$ is a rank 1
subsheaf of $V_P$, in which case (SB1)--(SB3) imply that
$\mathcal{S}'$ is not destabilizing.$\;\square$

\begin{Lem} Let
$(V,\theta)=V_P\stackrel{\oplus}{\longrightarrow}\mathcal{O}_X$ be
a stable projectively flat binary Higgs bundle such that ${\rm
im}(\theta|V_P)$ has rank $1$.  Let $V_1={\rm ker}(\theta|V_P)$
and $V_2={\rm im}(\theta|V_P)$. Then $V_1$ and $V_2$ are vertical
line bundles.

\label{bin_thm_2}
\end{Lem}
\emph{Proof.}$\;$We have an exact sequence $0 \rightarrow V_1
\rightarrow V_P \rightarrow V_2 \rightarrow 0.$  It follows that
there exist divisors $D_1$ and $D_2$ and a dimension 0 subscheme
$\tilde Z$ such that $V_1=\mathcal{O}_X(D_1)$ and $V_2=I_{\tilde
Z}\otimes\mathcal{O}_X(D_2)$, where $I_{\tilde Z}$ is the ideal
sheaf associated to $\tilde Z$ \cite{Friedman}.

We first show that $D_2$ is a vertical divisor.  Since $V_2$ is
the image of $\theta|V_P$, which maps to $\Omega^1_X$, we find
from the short exact sequence (\ref{short_exact_1}) in Lemma
\ref{tern_thm_1} that either ${\rm Hom}(I_{\tilde
Z}\otimes\mathcal{O}_X(D_2),N)\neq0$ or ${\rm Hom}(I_{\tilde
Z}\otimes\mathcal{O}_X(D_2),Q)\neq0$.  Since $\tilde Z$ has
codimension $2$, we then deduce that either
$H^0(\mathcal{O}_X(-D_2-2F+\sum(m_k-1)F_k))\neq 0$ or
$H^0(\mathcal{O}_X(-D_2+F))\neq 0$.  Let $H_0, H,$ and $k_0$ be as
in Def. \ref{Dolgachev_7}.  Since either $-D_2-2F+\sum(m_k-1)F_k$
or $-D_2+F$ is linearly equivalent to an effective divisor, we
have that $H_0\cdot D_2<\frac13 k_0<k_0$.  We also find that
$D_2\cdot F\leq 0$, with equality iff $D_2$ is vertical.  Suppose
that $D_2\cdot F<0$.  This would imply that $H\cdot D_2=H_0\cdot
D_2+k_0F\cdot D_2< k_0-k_0=0$.  But conditions (SB1)--(SB3) in
Lemma \ref{bin_thm_1} imply that $H\cdot D_2>0$.

We now show that $D_1$ is a vertical divisor.  We begin to do so
by showing that $F\cdot D_1=0$.  Suppose that $F\cdot D_1>0$.
Since $V$ is projectively flat, we have $(D_1+D_2)^2=3(\ell(\tilde
Z)+D_1\cdot D_2)$, where $\ell(\tilde Z)$ denotes the length of
$\tilde Z$ \cite{Friedman}. Since $D_2$ is vertical, we have that
$D_1\cdot D_2=\frac{(H_0\cdot D_2)(D_1\cdot F)}{H_0\cdot F}>0$.
It follows that $D_1^2>0$.  The Hodge index theorem, applied to
$(H\cdot F)D_1-(H\cdot D_1)F$, then shows that $H\cdot D_1\geq 0$.
But conditions (SB1)--(SB3) in Lemma \ref{bin_thm_1} imply that
$H\cdot D_1<0$.

Now suppose that $F\cdot D_1<0$.  This time, we apply the Hodge
index theorem to $(H_0\cdot F)D_1-(H_0\cdot D_1)F$ to find that
$H_0\cdot D_1\leq\frac{(H_0\cdot F)(3\ell(\tilde Z)+D_1\cdot
D_2)}{2(D_1\cdot F)}.$  It follows that $H\cdot D_1\leq
\frac{H\cdot D_2}{2}-k_0<-H\cdot D_2$.  But conditions
(SB1)--(SB3) in Lemma \ref{bin_thm_1} imply that $H\cdot
D_1>-H\cdot D_2$.

Hence $F\cdot D_1=0$.  The Hodge index theorem now implies that
$D_1^2\leq 0$ with equality iff $D_1$ is vertical.  Projective
flatness implies that $D_1^2=3\ell(\tilde Z)\geq 0$.  Therefore
$D_1$ is vertical.

Finally, since $D_1$ is vertical, we have that
$0=D_1^2=3\ell(\tilde Z)$, which implies that $V_2$ is a line
bundle.$\;\square$

\begin{Lem}
Let $V_1=\mathcal{O}_X(aF+\sum a_k F_k)$ and
$V_2=\mathcal{O}_X(bF+\sum b_k F_k)$ be vertical line bundles such
that $d_2\leq -2$, where $d_2=a-b-\#\{a_k<b_k\}$.

If there is a nonsplit extension of the
form\begin{equation}0\rightarrow V_1\rightarrow
V_P\stackrel{\beta}{\rightarrow} V_2 \rightarrow
0,\label{short_exact_5}\end{equation}and $L=\mathcal{O}_X(cF+\sum
c_k F_k)$ is a vertical line bundle with $d_1\geq 0$ and $d_3\leq
-2$, where $d_1=b-c-\#\{b_k<c_k\}$ and $d_3=a-c-\#\{a_k<c_k\}$,
such that $H^0(L^{^*}\!\otimes V_P)=0$, then
$d_1+1\leq\emph{min}(-d_2-1,-d_3-1)$.

Conversely, there exists a nonsplit extension {\rm
(\ref{short_exact_5})} such that if $L=\mathcal{O}_X(cF+\sum c_k
F_k)$ is any vertical line bundle with $d_1\geq 0$ and $d_3\leq
-2$ such that $d_1+1\leq\emph{min}(-d_2-1,-d_3-1)$, then
$H^0(L^{^*}\!\otimes V_P)=0$.

\label{bin_thm_18}
\end{Lem}
\emph{Proof.}  First, we show that if $V_P$ and $L$ are subject to
the given conditions, then $d_1+1\leq{\rm min}(-d_2-1,-d_3-1)$.

Observe that $L^{^*}\!\otimes V_1=\mathcal{O}_X(d_3F+\sum r_k
F_k)$ and $L^{^*}\!\otimes V_2=\mathcal{O}_X(d_1F+\sum r'_k F_k)$
for some $r_k, r'_k\geq 0$.  Consider the short exact sequence
\begin{equation}0\rightarrow L^{^*}\!\otimes V_1\rightarrow L^{^*}\!\otimes
V_P\rightarrow L^{^*}\!\otimes V_2\rightarrow
0.\label{short_exact_6}\end{equation}The associated long exact
sequence in cohomology then implies that the coboundary map
$\delta:H^0(L^{^*}\!\otimes V_2)\rightarrow H^1(L^{^*}\!\otimes
V_1)$ is injective.  Consequently, $h^0(L^{^*}\!\otimes V_2)\leq
h^1(L^{^*}\!\otimes V_1)$, and so by Lemma \ref{Dolgachev_2}, we
have that $d_1+1\leq -d_3-1$.

Suppose now that $d_1+1>-d_2-1$. Let $\sigma$ be an element of
$H^1(V_2^{^*}\!\otimes V_1)$ which defines the extension
(\ref{short_exact_5}). Note that $V_2^{^*}\!\otimes
V_1=\mathcal{O}_X(d_2F+\sum r_k F_k)$ for some $r_k$ with $r_k\geq
0$.  Taking notation from Lemma \ref{Dolgachev_6}(ii), we have
that $\sigma$ equals
$\sigma_{_{-1}}w_{\gamma}^{-1}+\dots+\sigma\!_{_{d_2+1}}w_{\gamma}^{d_2+1}$
on $V_{\gamma}\cap W_{\xi}$ and 0 elsewhere for some
$\sigma_{_{-1}},\dots,\sigma\!_{_{d_2+1}}$.  Let
$\{\phi'_{\alpha\beta}\}$ be a system of transition functions for
the line bundle $L^{^*}\!\otimes V_1$, and let
$\{\phi''_{\alpha\beta}\}$ be a system of transition functions for
the line bundle $L^{^*}\!\otimes V_2$.  We may regard $\sigma$ as
the extension class of (\ref{short_exact_6}).  Transition matrices
for $L^{^*}\!\otimes V_P$ are then given by $\left(
\begin{array}{cc}
\phi'_{\alpha\beta} & \phi''_{\alpha\beta}\sigma_{\alpha\beta} \\
0 & \phi''_{\alpha\beta} \\
\end{array} \right)$.

Let $s\in H^0(L^{^*}\!\otimes V_2)$ be the nonzero section such
that with respect to the trivialization on $V_{\gamma}$, we have
$s_{\gamma}=w_{\gamma}^{d_1}$, as in Lemma \ref{Dolgachev_6}(i).
Then
$\delta(s)=\sigma\!_{_{-1}}w_{\gamma}^{d_1-1}+\dots+\sigma\!_{_{d_2+1}}w_{\gamma}^{d_1+d_2+1}$
on $V_{\gamma}\cap W_{\xi}$ and 0 elsewhere.  So, by Lemma
\ref{Dolgachev_6}(ii) and the inequality $d_1+1>-d_2-1$, we have
that $\delta(s)=0\in H^1(L^{^*}\!\otimes V_1).$  But since
$\delta$ is injective, this yields the desired contradiction.

Conversely, we now show that there exists a nonsplit extension
(\ref{short_exact_5}) such that if $L=\mathcal{O}_X(cF+\sum c_k
F_k)$ is any vertical line bundle with $d_1\geq 0$ and $d_3\leq
-2$ such that $d_1+1\leq{\rm min}(-d_2-1,-d_3-1)$, then
$H^0(L^{^*}\!\otimes V_P)=0$.  Let
$(\sigma\!_{_{-1}},\sigma\!_{_{-2}},\dots,\sigma\!_{_{d_2+1}})$ be
a $(-d_2-2)$-tuple of complex numbers such that for any
$\ell_1,\ell_3$ with $\ell_1\geq 0$ and $\ell_3\leq -2$ such that
$\ell_1+1\leq{\rm min}(-d_2-1,-\ell_3-1)$, the
matrix$$\Theta_{\ell_1,\ell_3}=\left(
\begin{array}{cccccccc}
\sigma\!_{_{\ell_3+1}}     & \sigma\!_{_{\ell_3}} & \dots & \sigma\!_{_{d_2+1}} & 0                   & 0 & \dots & 0                    \\
\sigma\!_{_{\ell_3+2}}     & \sigma\!_{_{\ell_3+1}} & \dots & \sigma\!_{_{d_2+2}} & \sigma\!_{_{d_2+1}} & 0 & \dots & 0                    \\
\vdots                     & &       &                     &                     & \ddots  &       & \vdots               \\
\sigma\!_{_{\ell_3+d_2+1}} & \sigma\!_{_{\ell_3+d_2}} & \dots &                     & \dots               &   & \dots & \sigma\!_{_{d_2+1}}  \\
\vdots                     & &       &                     &               &   &       & \vdots               \\
\sigma\!_{_{-1}}           & \sigma\!_{_{-2}} & \dots &                     & \dots               &   & \dots & \sigma\!_{_{-\ell_1-1}} \\
\end{array} \right)$$has maximal rank.  (One may construct such a
sequence of $\sigma$'s by induction on $-d_2-1$; given
$\sigma\!_{_{-1}},\sigma\!_{_{-2}},\dots,\sigma\!_{_{d_2}}$,
choose $\sigma\!_{_{d_2+1}}$ so that every \emph{square} matrix of
the above form has nonzero determinant.  This is possible because
there are only finitely many such matrices, and for each such
matrix, the determinant is zero for only finitely many values of
$\sigma\!_{_{d_2+1}}$.)

Let $\sigma$ be the element in $H^1(V_2^{^*}\!\otimes V_1)$
represented by a 1-cocycle which equals
$\sigma_{\gamma\xi}=\sigma_{_{-1}}w_{\gamma}^{-1}+\dots+\sigma\!_{_{d_2+1}}w_{\gamma}^{d_2+1}$
on $V_{\gamma}\cap W_{\xi}$ and 0 elsewhere.  Let $V_P$ be the
rank 2 bundle given as an extension as in (\ref{short_exact_5})
whose extension class is determined by $\sigma$.  Since $\sigma$
is nonzero, (\ref{short_exact_5}) does not split.  Let
$L=\mathcal{O}_X(cF+\sum c_k F_k)$ be a vertical line bundle with
$d_1\geq 0$ and $d_3\leq -2$ such that $d_1+1\leq{\rm
min}(-d_2-1,-d_3-1)$.  We must show that $H^0(L^{^*}\!\otimes
V_P)=0$.

The condition $d_3\leq -2$ guarantees that $H^0(L^{^*}\!\otimes
V_1)=0$.  It therefore suffices to show that the coboundary map
$\delta:H^0(L^{^*}\!\otimes V_2)\rightarrow H^1(L^{^*}\!\otimes
V_1)$ is injective.  Let $s\in H^0(L^{^*}\!\otimes V_2)$.  We now
show that if $\delta(s)=0$, then $s=0$.

From Lemma \ref{Dolgachev_6}(i), we know that on $V_{\gamma}$, the
section $s$ is of the form $s_{\gamma}=s_0+s_1
w_{\gamma}+\dots+s_{d_1} w_{\gamma}^{d_1}$ with respect to the
trivialization on $V_{\gamma}$.  From Lemma \ref{Dolgachev_6}(ii),
we know that if $c$ is the 1-cocycle given by $w^j$ on
$V_{\gamma}\cap W_{\xi}$ and 0 elsewhere, then $[c]=0\in
H^1(L^{^*}\!\otimes V_1)$ if and only if $j\geq 0$ or $j\leq
-d_3$. Since $\delta(s)$ equals $s_{\gamma}\sigma_{\gamma\xi}$ on
$V_{\gamma}\cap W_{\xi}$ and 0 elsewhere, we have that
$\delta(s)=0$ if and only if the following equalities
hold:$$\sigma\!_{_{d_3+1}}s_{_{0}}+\sigma\!_{_{d_3}}s_{_{1}}+\dots+\sigma_{_{d_2+1}}s_{_{d_3-d_2}}=0$$
$$\sigma\!_{_{d_3+2}}s_{_{0}}+\sigma\!_{_{d_3+1}}s_{_{1}}+\dots+\sigma_{_{d_2+1}}s_{_{d_3-d_2+1}}=0$$
$$\dots$$
$$\sigma\!_{_{d_3+d_2+1}}s_{_{0}}+\sigma\!_{_{d_3+d_2}}s_{_{1}}+\dots+\sigma_{_{d_2+1}}s_{_{d_1}}=0$$
$$\dots$$
$$\sigma\!_{_{-1}}s_{_{0}}+\sigma\!_{_{-2}}s_{_{1}}+\dots+\sigma_{_{-d_1-1}}s_{_{d_1}}=0$$

Since $\Theta_{d_1,d_3}$ has maximal rank and $d_1+1\leq -d_3-1$
(which is to say, regarding the $s$'s as variables, that there are
at least as many equations as variables), we conclude that
$s=0$.$\;\square$

\emph{Proof of Prop. \!\ref{bin_thm_20}.} We first show that if
$(V,\theta)=V_P\stackrel{\oplus}{\longrightarrow}\mathcal{O}_X$ is
a stable projectively flat binary Higgs bundle with
$\rm{rank}(\rm{im}(\theta|V_P))=1$, then $V_P$ has the stated
form.

Lemma \ref{bin_thm_2} implies that $V_1=\rm{ker}(\theta|V_P)$ and
$V_2=\rm{im}(\theta|V_P)$ are vertical line bundles; we therefore
obtain the extension (\ref{short_exact_4}).  Condition (i) follows
from (SB1) and (SB3) of Lemma \ref{bin_thm_1}.

Stability implies that (\ref{short_exact_4}) does not split;
therefore $h^1(V_2^{^*}\!\otimes V_1)>0$. It follows from (i) that
$d_2<0$.  Condition (ii) then follows from Lemma
\ref{Dolgachev_2}(ii).

Since $V_2$ is a subsheaf of $\Omega^1_X$, we must have that
$H^0(V_2^{^*}\!\otimes\Omega^1_X)\neq 0$.  Condition (iii) then
follows from Lemma \ref{tern_thm_7}.

Let $(c,c_1,\dots,c_n)$ be an $(n+1)$-tuple of integers such that
$0\leq c_k<m_k$ for all $k$ and $d_1\geq 0$ and
$C\geq\frac{2}{3}(A+B)$.  Let $L=\mathcal{O}_X(cF+\sum c_k F_k)$.
From (SB2) of Lemma \ref{bin_thm_1}, we know that
$H^0(L^{^*}\!\otimes V_P)=0$.  Note that $L^{^*}\!\otimes
V_1=\mathcal{O}_X(d_3F+\sum r_k F_k)$ for some $r_k$ with $0\leq
r_k<m_k$ for all $k$.  Arguing as in the proof that condition (ii)
holds, we see that $d_3<0$. From the long exact sequence in
cohomology associated to (\ref{short_exact_6}), we find that
$H^1(L^{^*}\!\otimes V_1)\neq 0$. Lemma \ref{Dolgachev_2} then
implies that $d_3\leq -2$.  Condition (iv) then follows from Lemma
\ref{bin_thm_18}.

Conversely, suppose that we are given $a$'s and $b$'s satisfying
conditions (i)--(iv), and let $V_1=\mathcal{O}_X(aF+\sum a_k F_k)$
and $V_2=\mathcal{O}_X(bF+\sum b_k F_k)$.  We will show that there
exists a stable projectively flat binary Higgs bundle
$V_P\stackrel{\oplus}{\longrightarrow}\mathcal{O}_X$ with
$\rm{rank}(\rm{im}(\theta|V_P))=1$ and $V_P$ as in
(\ref{short_exact_4}).

Lemma \ref{bin_thm_18} and condition (ii) guarantee the existence
of a rank 2 bundle $V_P$ and a nonsplit extension
(\ref{short_exact_4}) such that if $L=\mathcal{O}_X(cF+\sum c_k
F_k)$ is any vertical line bundle with $d_1\geq 0$ and $d_3<0$
such that $d_1+1\leq{\rm min}(-d_2-1,-d_3-1)$, then
$H^0(L^{^*}\!\otimes V_P)=0$. By Lemma \ref{tern_thm_7} and
condition (iii), there exists a nonzero map
$\alpha:V_2\rightarrow\Omega^1_X$.  Let
$V=V_P\oplus\mathcal{O}_X$.  Define a Higgs field $\theta$ by
$\theta|V_P=\alpha\circ\beta$ and $\theta|\mathcal{O}_X=0$. Note
that $\theta\wedge\theta=0$. Then $(V,\theta)$ is a binary Higgs
bundle with $\rm{rank}(\rm{im}(\theta|V_P))=1$.  Moreover, $V$ is
projectively flat, since $0=c_1^2(V)=3c_2(V).$

It remains to be shown that $(V,\theta)$ is stable.  (SB1) and
(SB3) from Lemma \ref{bin_thm_1} follow from condition (i).  Let
us now verify that (SB2) holds.  Suppose to the contrary that
there exists a rank 1 subsheaf $\mathcal{S}$ of $V_P$ such that
${\rm deg}(\mathcal{S})\geq\frac{2}{3}\,{\rm deg}(V_P)$.  Let $L$
be the kernel of the natural map $V_P\rightarrow
\frac{\frac{V_P}{\mathcal{S}}}{{\rm
Tor}(\frac{V_P}{\mathcal{S}})}$. Then $L$ is a line bundle, ${\rm
deg}(L)\geq\rm{deg}(\mathcal{S})$, and $H^0(L^{^*}\!\otimes
V_P)\neq 0$.  (See \cite{Kobayashi}.)  Stability, together with
Def. \ref{Dolgachev_7}, implies that $L$ is vertical.

Write $L=\mathcal{O}_X(cF+\sum c_k F_k)$.  Dividing both sides of
${\rm deg}(L)\geq\frac{2}{3}\,{\rm deg}(V_P)$ by $H\cdot F$, we
find that $C\geq\frac{2}{3}(A+B)$.  Note that $L^{^*}\!\otimes
V_2=\mathcal{O}_X(d_1F+\sum r_k F_k)$, and so $H^0(L^{^*}\!\otimes
V_2)\neq 0$ implies that $d_1\geq 0$.  It now follows from
condition (iv) that $d_1+1\leq{\rm min}(-d_2-1,-d_3-1)$. Moreover,
$d_3<0$ since $H^0(L^{^*}\!\otimes V_1)=0$.  Our choice of $V_P$
then implies that $H^0(L^{^*}\!\otimes V_P)=0$, contradicting our
earlier assertion that $H^0(L^{^*}\!\otimes V_P)\neq 0$. Therefore
$(V,\theta)=V_P\stackrel{\oplus}{\longrightarrow}\mathcal{O}_X$ is
stable, as desired.$\;\square$

\subsection{The case of the stable binary Higgs bundle with rank(im($\theta$))=2}
\label{stable_binary_rank_2}

~\newline In this subsection, we show that there does not exist a
stable binary Higgs bundle $(V,\theta)$ on $X$ with
rank(im($\theta$))=2. Throughout this section, let
$N=\mathcal{O}_X (-2F+\sum_k (m_k-1)F_k)$ and
$Q=I_Z\otimes\mathcal{O}_X(F)$, as in Lemma \ref{tern_thm_1}.

\begin{Lem}
Suppose that
$(V,\theta)=V_P\stackrel{\oplus}{\longrightarrow}\mathcal{O}_X$ is
a stable projectively flat binary Higgs bundle with
$\emph{rank}(\emph{im}(\theta))=2$.  Then there exists an exact
sequence $$0\to V_1\to V_P\to V_2\to 0,$$ where $V_1$ and $V_2$
are vertical line bundles and $H^0(V^{^*}_2\otimes Q)\neq 0$.

\label{bin_thm_25}
\end{Lem}
\emph{Proof.} Let $\beta$ be the map in the exact sequence of
Lemma \ref{tern_thm_1} from $\Omega^1_X$ to $Q$.  Let $V_2={\rm
im}(\beta\circ(\theta|V_P))$, and let $V_1={\rm
ker}(\beta\circ(\theta|V_P))$.  This gives us an exact sequence
$0\to V_1\to V_P\to V_2\to 0$.  Since
$\rm{rank}(\rm{im}(\theta))=2$, we see that
$1=\rm{rank}(V_2)=\rm{rank}(V_1)$.  The proof of Lemma
\ref{bin_thm_2} shows that $V_1$ and $V_2$ are vertical line
bundles.  Moreover, the inclusion map $\iota:V_2\hookrightarrow Q$
yields a nonzero element of $H^0(V^{^*}_2\otimes Q)$.$\;\square$

\begin{Prop}
If $(V,\theta)$ is a stable binary Higgs bundle, then
$\rm{im}(\theta)$ has rank $1$.

\label{bin_thm_10}
\end{Prop}
\emph{Proof.}  By tensoring with a line bundle, as in Lemma
\ref{Higgs_5}, we may assume that $(V,\theta)$ is of the form
$V_P\stackrel{\oplus}{\longrightarrow}\mathcal{O}_X$.  Then
$\rm{im}(\theta)$ is a subsheaf of $\Omega^1_X$ and so has rank 0,
1, or 2.  As noted in the introduction to
\S\ref{stable_binary_rank_1}, $\rm{im}(\theta)$ cannot have rank
$0$.

Suppose $\rm{im}(\theta)$ has rank 2.  By Lemma \ref{bin_thm_25},
we have an exact sequence $$0\to V_1\to V_P\to V_2\to 0,$$ where
$V_1$ and $V_2$ are vertical line bundles and $H^0(V^{^*}_2\otimes
Q)\neq 0$.  By Lemma \ref{tern_thm_2}, we have that ${\rm
deg}(V_2)<0$. As in Lemma \ref{bin_thm_1}, stability implies that
${\rm deg}(V_P)>0$, whence we see that $0<{\rm deg}(V_P)={\rm
deg}(V_1)+{\rm deg}(V_2)<{\rm deg}(V_1).$ The proof of Lemma
\ref{bin_thm_1} also shows that ${\rm deg}(V_1)<\frac{2}{3}\,{\rm
deg}(V_P)$, whereby one obtains the contradictory inequality
$$0<{\rm deg}(V_1)<2\,{\rm deg}(V_2)<0.\;\;\;\square$$

\subsection{The case of the reducible ternary Higgs bundle}
\label{red_tern}

~\newline We now consider reducible, polystable, ternary Higgs
bundles of the form
$(V,\theta)=V_2\stackrel{\oplus}{\rightarrow}V_3
\stackrel{\oplus}{\rightarrow} V_1$.  In this case, either
$\theta|V_2$ or $\theta|V_3$ must be the zero map.  (For if not,
then $V$ is not reducible.)  We divide into three cases
accordingly, depending whether the first map only is zero, the
second map only is zero, or both
are.\vspace{-.1in}\newline\newline\textbf{Case 1}: $\theta|V_2=0$
and $\theta|V_3\neq 0$

\begin{Prop} There exists
a polystable ternary Higgs bundle
$(V,\theta)=V_2\stackrel{\oplus}{\rightarrow}V_3
\stackrel{\oplus}{\rightarrow} V_1$ with $\theta|V_2=0$ and
$\theta|V_3\neq 0$ and $c_1(V_2)=c_1(V_3\oplus V_1)=c_2(V_3\oplus
V_1)=0$ if and only if $V_2=\mathcal{O}_X$ and
$V_3=\mathcal{O}_X(bF+\sum b_k F_k)$ and $V_1=V_3^{^*}\!$, where
the $b$'s are subject to the following numerical conditions:

{\rm (i)} $\;B=b+\sum\frac{b_k}{m_k}>0$, and

{\rm (ii)} $\;2b+\#\{b_k\geq\frac{m_k}{2}\}\leq -2.$

\label{red_tern_10}
\end{Prop}
\emph{Proof.}  First, let $V_2=\mathcal{O}_X$ and
$V_3=\mathcal{O}_X(bF+\sum b_k F_k)$ and $V_1=V_3^{^*}\!$, where
the $b$'s satisfy (i) and (ii).  Note that $V_3\otimes
V_3=\mathcal{O}_X((2b+\#\{b_k\geq\frac{m_k}{2}\})F+\sum r_k F_k)$
for some $r_k$ with $0\leq r_k<m_k$.  Condition (ii) guarantees
that there exists a nonzero map $\theta:V_3\rightarrow
V_1\otimes\Omega^1_X$, by Lemma \ref{tern_thm_7}. Extend $\theta$
to $V$ by letting $\theta|V_2=\theta|V_1=0$; then
$\theta\wedge\theta=0.$ Condition (i) guarantees that
$V_3\stackrel{\oplus}{\longrightarrow}V_1$ is stable. We have
$c_1(V_3\oplus V_1)=0$ since $V_1=V_3^{^*}\!$.  Also,
$c_2(V_3\oplus V_1)=c_1(V_3)c_1(V_1)=0$ since $V_3$ and $V_1$ are
vertical.

Now let $(V,\theta)=V_2\stackrel{\oplus}{\rightarrow}V_3
\stackrel{\oplus}{\rightarrow} V_1$ be a polystable ternary Higgs
bundle with $\theta|V_2=0$ and $\theta|V_3\neq 0$ and
$c_1(V_2)=c_1(V_3\oplus V_1)=c_2(V_3\oplus V_1)=0$.  Then
$V_2=\mathcal{O}_X$ and $V_1=V_3^{^*}\!$, since $c_1(V_2)=0$ and
$c_1(V_3\oplus V_1)=0$.  Write $V_3=\mathcal{O}_X(D)$ for some
divisor $D$.  We now show that $D$ is linearly equivalent to
$bF+\sum b_k F_k$ for some $n+1$-tuple of $b$'s.  (In other words,
we show that $D$ is vertical.)  First, we show that $D\cdot F=0$.
 Suppose that $D\cdot F>0$.  Let $H_0, H,$ and $k_0$ be as in Def.
\ref{Dolgachev_7}. Using the same line of reasoning as in the
proof of Lemma \ref{bin_thm_1}, we see that the nonvanishing of
$\theta|V_3$ implies that $H\cdot D<k_0$.  From this we find that
$H_0\cdot D<0$.  The condition $c_2(V_3\oplus V_1)=0$ implies that
$D^2=0$.  The Hodge index theorem, applied to $(H_0\cdot
F)D-(H_0\cdot D)F$, then yields the contradictory inequality
$-2(H_0\cdot F)(H_0\cdot D)(D\cdot F)\leq 0$.  Now suppose that
$D\cdot F<0$.  This time, we use the fact that the nonvanishing of
$\theta|V_3$ implies that $H_0\cdot D<k_0$, which shows in turn
that $H\cdot D=H_0\cdot D+k_0 F\cdot D<0$. But stability implies
the contradictory inequality $H\cdot D>0$. Therefore $D\cdot F=0$.
The Hodge index theorem, applied to $(H\cdot F)D-(H\cdot D)F$,
then implies that $D$ is vertical. We obtain condition (i) by
diving both sides of the inequality ${\rm deg}(V_3)>0$ by $H\cdot
F$. Lemma \ref{tern_thm_7} yields condition
(ii).$\;\square$\vspace{-.1in}\newline\newline\textbf{Case 2}:
$\theta|V_2\neq 0$ and $\theta|V_3=0$.

This case is the same as Case 1, with the $V$'s
relabeled.\vspace{-.1in}\newline\newline\textbf{Case 3}:
$\theta|V_2=\theta|V_3=0$

This case is trivial; there exists a polystable Higgs bundle
$V_2\stackrel{\oplus}{\rightarrow}V_3
\stackrel{\oplus}{\rightarrow} V_1$ with
$c_1(V_2)=c_1(V_3)=c_1(V_1)=0$ and $\theta|V_2=\theta|V_3=0$ if
and only if $V_2=V_3=V_1=\mathcal{O}_X$.

\subsection{The case of the reducible binary Higgs bundle}
\label{red_bin}

Let $(V,\theta)=V_P\stackrel{\oplus}{\longrightarrow}V_Q$ be a
reducible polystable binary Higgs bundle whose summands have
vanishing Chern classes, where rank($V_P$)=2 and rank($V_Q$)=1.
The rank $R$ of the image of $\theta$ in $V_Q\otimes\Omega^1_X$ is
either 2, 1, or 0. If $R=2$, then $(V,\theta)$ can not be
reducible.  If $R=1$, then we must have $V_P=V_1\oplus V_2$, where
$V_1={\rm ker}(\theta|V_P)$; this case was discussed in
\S\ref{red_tern}.  If $R=0$, then $\theta$ is the zero map.  In
this case, we must have $V_Q=\mathcal{O}_X$ and $V_P$ stable.

\textbf{Remark:} An explicit description of all stable rank 2
bundles on $X$ with vanishing Chern classes can be found in
\cite[Proposition 4.1]{BO}.  (The method of proof of Prop.
\ref{bin_thm_20} also yields such a description.)

\section{Main Theorem and an Example}
\label{main}

Putting together the pieces from the previous sections, we have
the following explicit description of all orbifold Toledo
invariants that arise from semisimple ${\rm U}(2,1)$
representations of the orbifold fundamental group of the
2-orbifold associated to a Seifert fibered homology 3-sphere.

\begin{Thm}
Let $O$ be the base orbifold of a large \cite[\S 5.3]{Orlik}
Seifert fibered homology $3$-sphere.  Let $n$ equal the number of
cone points that $O$ has, and let $m_1,\dots,m_n$ denote the
orders of these cone points. Then there exists a semisimple
representation $\rho:\pi_1^{\rm{orb}}(O)\rightarrow \rm{U(2,1)}$
such that $\tau=\tau_{\rm{orb}}(\rho)$ if and only if
$\tau=\pm(a+b+\sum\frac{a_k+b_k}{m_k})$ for some $(2n+2)$-tuple
$(a, a_1, \dots, a_n, b, b_1, \dots, b_n)$ of integers with $0\leq
a_k,b_k<m_k$ for all $k=1,\dots,n$ such that at least one of
\emph{(i)--(iv)} holds:

\emph{(i)} The $a$'s and $b$'s satisfy {\rm (i)--(iv)} from Prop.
{\rm \ref{tern_thm_11}} as well as \emph{($\star$)} below; or

\emph{(ii)} The $a$'s and $b$'s satisfy {\rm (i)--(iv)} from Prop.
{\rm \ref{bin_thm_20}} as well as \emph{($\star$)} below; or

\emph{(iii)}  The $b$'s satisfy {\rm (i)} and {\rm (ii)} from
Prop. {\rm \ref{red_tern_10}}, and $a=a_k=0$ for all $k$; or

\emph{(iv)} $a=b=a_k=b_k=0$ for all $k$.

\emph{($\star$)} There exist integers
$y,y_1,\dots,y_n,s_1,\dots,s_n$ such that $3y+(s_1+\dots+s_n)=a+b$
and $3y_k-m_k s_k=a_k+b_k$ for $k=1,\dots,n$.

\label{main_1}
\end{Thm}
\emph{Proof.} By Lemmas \ref{Higgs_3}, \ref{Higgs_5}(ii),
\ref{Hodge_2}, and \ref{orbifold_7}, as well as the discussion
following Def. \ref{Hodge_3}, it suffices to show that
$c_1(\mathcal{O}_X((a+b)F+\sum(a_k+b_k)F_k))$ equals the Higgs
bundle Toledo invariant of a stable ternary, stable binary,
reducible ternary, or reducible binary Higgs bundle whose summands
have vanishing Chern classes if and only if the $a$'s and $b$'s
satisfy one of (i)--(iv).

Suppose $(V,\theta)=V_2\stackrel{\oplus}{\rightarrow}V_3
\stackrel{\oplus}{\rightarrow} V_1$ is a stable ternary Higgs
bundle with vanishing Chern classes.  By Lemma \ref{Higgs_5},
tensoring with $V_3^{^*}\!$ yields a stable ternary Higgs bundle
$(V',\,\theta')=(V\otimes V_3^{^*}\!,\,\theta\otimes
1)=(V_2\otimes
V_3^{^*}\!)\stackrel{\oplus}{\rightarrow}\mathcal{O}_X
\stackrel{\oplus}{\rightarrow} (V_2\otimes V_3^{^*}\!)$ with
$\tau_{(V,\,\theta)}=\tau_{(V',\,\theta')}$. By Prop.
\ref{tern_thm_11} and Definition \ref{Higgs_4}, we then have that
$\tau_{(V',\,\theta')}=c_1(\mathcal{O}_X((a+b)F+\sum(a_k+b_k)F_k))$,
where the $a$'s and $b$'s satisfy (i)--(iv) from Prop.
\ref{tern_thm_11}.  Moreover,
$\mathcal{O}_X((a+b)F+\sum(a_k+b_k)F_k)={\rm
det}(V')=V_3^{^*}\!\otimes V_3^{^*}\!\otimes V_3^{^*}\!$ is
vertical.  Thus $V_3^{^*}\!$ is of the form $\mathcal{O}_X(yF+\sum
y_k F_k)$ with $3(yF+\sum y_k F_k)$ linearly equivalent to
$(a+b)F+\sum(a_k+b_k)F_k$. Condition ($\star$), which is
equivalent to the condition that $(a+b)F+\sum(a_k+b_k)F_k$ is
``divisible by 3,'' therefore holds.

Conversely, given $a$'s and $b$'s satisfying (i), Prop.
\ref{tern_thm_11} and Def. \ref{Higgs_4} guarantee the existence
of a stable projectively flat ternary Higgs bundle $(V',\theta')$
with
$\tau_{(V',\,\theta')}=\mathcal{O}_X((a+b)F+\sum(a_k+b_k)F_k)$.
Condition ($\star$) is then equivalent to the existence of a
vertical line bundle $V_3=\mathcal{O}_X(yF+\sum y_k F_k)$ such
that $c_1(V'\otimes V_3)=c_2(V'\otimes V_3)=0$.  By Lemma
\ref{Higgs_5}, $V'\otimes V_3$ is a stable ternary Higgs bundle
with $\tau_{(V,\theta)}=\tau_{(V',\,\theta')}$.

To summarize: $c_1(\mathcal{O}_X((a+b)F+\sum(a_k+b_k)F_k))$ equals
the Higgs bundle Toledo invariant of a stable flat ternary Higgs
bundle on $X$ if and only if the $a$'s and $b$'s satisfy (i).

A similar argument, using Prop. \ref{bin_thm_20} instead of Prop.
\ref{tern_thm_11}, shows that
$c_1(\mathcal{O}_X((a+b)F+\sum(a_k+b_k)F_k))$ equals the Higgs
bundle Toledo invariant of a stable binary Higgs bundle
$(V,\theta)$ with $c_1(V)=c_2(V)=0$ and rank(im($\theta$))=1 if
and only if the $a$'s and $b$'s satisfy (ii).  Prop.
\ref{bin_thm_10} shows that there are no stable binary Higgs
bundles with rank(im($\theta$))=2.  Condition (iii) covers Cases 1
and 2 from \S\ref{red_tern}, and (iv) covers Case 3 from
\S\ref{red_tern} as well as the reducible binary case with
$\theta=0$ (as discussed in \S\ref{red_bin}), since the Toledo
invariant vanishes in both of these cases.$\;\square$

\begin{Cor}

{\rm (a)} A lower bound for the number of distinct connected
components in the representation space
$\mathcal{R}^+_{\rm{U}(2,1)}(O)=\frac{\rm{Hom}^+(\pi^{\rm{orb}}_1(\it{O}),\,\rm{U}(2,1))}{\rm{U}(2,1)}$
is given by the number of distinct values
$\,\pm(a+b+\sum\frac{a_k+b_k}{m_k})$, where the $a$'s and $b$'s
satisfy one of {\rm (i)--(iv)} from Thm. {\rm \ref{main_1}}.

{\rm (b)} A lower bound for the number of distinct connected
components in the representation space
$\mathcal{R}^{^*}\!\!_{\rm{PU}(2,1)}(Y)=\frac{{\rm
Hom}^*(\pi_1(\it{Y}),\,\rm{PU}(2,1))}{\rm{PU}(2,1)}$ is given by
the number of distinct values
$\,\pm(a+b+\sum\frac{a_k+b_k}{m_k})$, where the $a$'s and $b$'s
satisfy {\rm (i)} or {\rm (ii)} from Thm. {\rm \ref{main_1}}.

\label{main_2}
\end{Cor}
\emph{Proof.} We prove (b) only; the proof of (a) is similar.
Lemma \ref{Seifert_2} shows that we may replace $Y$ by $X$ in the
statement of this theorem. Lemma \ref{Toledo_2} shows that
(equivalence classes of) PU(2,1) representations with distinct
Toledo invariants lie in distinct components of
$\mathcal{R}^{^*}\!\!_{\rm{PU}(2,1)}(X)$. If $\rho\in{\rm
Hom}^*(\pi_1(\it{X}),\,\rm{U}(2,1))$, then
$\varphi\circ\rho\in{\rm Hom}^*(\pi_1(\it{X}),\,\rm{PU}(2,1))$,
where $\varphi:\rm{U}(2,1)\rightarrow\rm{PU}(2,1)$ is the
canonical homomorphism.  Lemmas \ref{Higgs_2} and \ref{Higgs_3}
show that the number of distinct Toledo invariants arising from
irreducible U(2,1) representations of $\pi_1(X)$ exactly equals
the number of distinct Higgs bundle Toledo invariants of stable
U(2,1) system of Hodge bundles on $X$ with vanishing Chern
classes. There exist $a$'s and $b$'s satisfying (i) or (ii) from
Thm. \ref{main_1} if and only if $\pm
c_1(\mathcal{O}_X((a+b)F+\sum(a_k+b_k)F_k))$ equals the Higgs
bundle Toledo invariant of a stable U(2,1) Higgs bundle on $X$
with vanishing Chern classes---in which case, by Lemma
\ref{orbifold_7}, the corresponding orbifold Toledo invariant is
$\pm(a+b+\sum\frac{a_k+b_k}{m_k})$.$\;\square$

\textbf{Example:} Let $n=3$, and let $(m_1,m_2,m_3)=(2,3,11)$.
Departing from our previous notations, let $F_{m_k}$ (instead of
$F_k$) denote the multiple fibre on $X$ of multiplicity $m_k$.

Let
$(V_1,\theta_1)=\mathcal{O}_X(-2F+F_2+2F_3+10F_{11})\stackrel{\oplus}{\rightarrow}
\mathcal{O}_X\stackrel{\oplus}{\rightarrow}\mathcal{O}_X(-F+F_2+F_3+F_{11})$
be a stable ternary Higgs bundle.

Let $(V_2,\theta_2)=\mathcal{O}_X\stackrel{\oplus}{\rightarrow}
\mathcal{O}_X\stackrel{\oplus}{\rightarrow}\mathcal{O}_X$, where
$\theta_2$ is the zero map.

Let $(V_3,\theta_3)$ be a stable binary Higgs bundle of the form
$V_P\stackrel{\oplus}{\rightarrow}\mathcal{O}_X$, where $V_P$ is
given by a nontrivial extension
$$0\to\mathcal{O}_X(-F+F_3+7F_{11})\to
V_P\to\mathcal{O}_X(-2F+F_2+2F_3+10F_{11})\to 0.$$

Let
$(V_4,\theta_4)=\mathcal{O}_X(-2F+F_2+2F_3+10F_{11})\stackrel{\oplus}{\rightarrow}
\mathcal{O}_X\stackrel{\oplus}{\rightarrow}\mathcal{O}_X(-F+F_2+F_3+2F_{11})$
be a stable ternary Higgs bundle.

Theorem \ref{main_1} guarantees that all orbifold Toledo
invariants arise from these four Higgs bundles and their duals.
Let $\tau_k$ be the orbifold Toledo invariant corresponding to
$(V_k,\theta_k)$.  Then $0=\tau_1=\tau_2$, $0.0455\approx\tau_3$,
and $0.0909\approx\tau_4$.  We conclude that in this case,
$\mathcal{R}^+_{\rm{U}(2,1)}(O)$ contains at least $5$ distinct
connected components.

Though $(V_1,\theta_1)$ and $(V_2,\theta_2)$ have the same Higgs
bundle Toledo invariant, it is unclear whether they lie in the
same component of $\mathcal{M}_{\rm{Dol}}(\rm{U}(2,1))$.  We hope
to address this question in a future paper.

\bibliographystyle{plain}
\bibliography{Bib}

\end{document}